\documentclass[11pt]{article}
\usepackage{mathtools}
\usepackage[colorlinks=true,linkcolor=blue,citecolor=blue,urlcolor=blue]{hyperref}
\usepackage[utf8]{inputenc}
\usepackage{amssymb,amsmath,amsfonts,eucal,mathrsfs,amsthm}
\usepackage[colorinlistoftodos]{todonotes}
\usepackage{xcolor}
\usepackage{enumitem}
\usepackage{accents}
\newtheorem{theorem}{Theorem}
\newtheorem{proposition}[theorem]{Proposition}
\newtheorem{lemma}[theorem]{Lemma}

\newtheorem{corollary}[theorem]{Corollary}

\newtheorem{remark}[theorem]{Remark}

\newtheorem{example}[theorem]{Example}

\theoremstyle{definition}
\newcommand{\R}{\mathbb{R}}
\newcommand{\Z}{\mathbb{Z}}
\newcommand{\Q}{\mathbb{Q}}
\newcommand{\Sf}{\mathbb{S}}
\newcommand{\C}{\mathbb{C}}

\newcommand{\spa}{\mbox{span}}

\newcommand{\Ric}{\mbox{Ric}}

\def\<{{\langle}}
\def\>{{\rangle}}
\def\B{\mathcal{B}}

\def\CP{\mathord{\mathbb C}\mathord{P}}

\def\n{\nabla}

\def\a{\alpha}

\def\be{\begin{equation}}
\def\ee{\end{equation}}

\begin{document}

\title{Rigidity of submanifolds under extrinsic curvature pinching}

\author{Theodoros Vlachos}
\date{}
\maketitle

\renewcommand{\thefootnote}{\fnsymbol{footnote}} 
\renewcommand{\thefootnote}{\arabic{footnote}} 

\renewcommand{\thefootnote}{\fnsymbol{footnote}} 
\renewcommand{\thefootnote}{\arabic{footnote}}

\begin{abstract}

We investigate the rigidity of compact submanifolds 
of arbitrary codimension in Riemannian manifolds 
whose curvature operator is bounded below by a 
nonnegative constant, under pointwise pinching 
conditions involving the second fundamental form 
and the mean curvature. In contrast to many earlier 
results, we allow equality in these conditions without 
imposing additional assumptions on the mean curvature 
or restricting the fundamental group of the 
submanifold. Allowing equality reveals geometric 
and topological phenomena excluded by strict pinching. 
We establish classification 
theorems in all dimensions $n\geq4$, describing 
the possible intrinsic geometries and topologies 
of such submanifolds. Our results also provide 
substantial information about the extrinsic geometry 
of the immersions, including precise conclusions 
concerning relative nullity, Dupin principal normals, 
and total geodesicity. As applications, we derive 
rigidity and sphere theorems for submanifolds 
of space forms and construct families of examples 
showing that the pinching conditions do not, in 
general, imply extrinsic rigidity. These results 
extend and refine several known rigidity theorems 
by providing a classification of the equality 
cases in a general ambient setting. 
\end{abstract}

\tableofcontents

\section{Introduction}
Curvature pinching inequalities often impose strong 
restrictions on the geometry and topology of Riemannian 
manifolds. For isometrically immersed submanifolds, 
pinching inequalities involving the second fundamental 
form and the mean curvature relate these intrinsic 
restrictions to the extrinsic geometry of the immersion. 
A central question is to what extent rigidity persists 
when equality is allowed in the pinching condition.

The subject has its origins in the work of Simons 
\cite{Sim} on minimal submanifolds of spheres and in 
the rigidity theorem of Chern, do Carmo, and Kobayashi 
\cite{CdCK}. Subsequent developments have focused 
primarily on submanifolds of space forms; see, for 
instance, 
\cite{adc,cuisun,font,HW,a,OV,p,SX,v3,v4,Xu,XGu,XG} 
and the references therein. The equality case is 
particularly subtle, as it may give rise to geometric 
and topological configurations excluded by strict 
pinching.

The aim of this paper is to investigate the rigidity 
of isometric immersions 
$f\colon M^n\to\bar M^{n+m}$ of compact Riemannian 
manifolds of dimension $n\geq4$ and arbitrary codimension 
under extrinsic curvature pinching. We assume that the 
curvature operator of the ambient manifold is bounded 
below by a nonnegative constant $c$, without requiring 
constant sectional curvature or local symmetry. Our 
principal pinching condition is
\begin{equation}\label{1}
S\leq 4c+\frac{n^2H^2}{n-2},
\tag{$*$}
\end{equation}
where $S=\|\alpha_f\|^2$ is the squared norm of the 
second fundamental form 
$\alpha_f\colon TM\times TM\to N_fM$, and 
$H=\|\mathcal H_f\|$ is the length of the mean 
curvature vector 
$\mathcal H_f=\frac1n\operatorname{tr}\alpha_f$. 
Observe that the pinching condition is invariant under 
homothetic rescaling of the metrics, with the 
corresponding rescaling of $c$.
\flushbottom

In the space-form setting, Howard and Wei \cite{HW} 
studied the strict pinching inequalities
\[
S<2kc+\frac{n^2H^2}{n-k},\qquad 1\leq k\leq n/2,
\]
where $k$ is an integer, and obtained topological 
restrictions through the nonexistence of stable currents 
or stable submanifolds of dimension $k$. For $k=2$, 
their condition is precisely the strict counterpart 
of \eqref{1}. Li, Xu, and Zhao \cite{LXZ} established 
sphere theorems under \eqref{1} for submanifolds of 
space forms subject to additional hypotheses, while 
Betti-number estimates for $c=0$ were obtained in 
\cite{OV} using Bochner methods. Related pinching 
problems in complex and quaternionic projective 
spaces, as well as in K\"ahler ambient manifolds, 
were investigated in \cite{ps,ss1,ss2}.

The distinction between strict and nonstrict pinching 
is central to our results. Condition~\eqref{1} implies 
that $M^n$ has nonnegative isotropic curvature, but 
equality allows the isotropic curvature to vanish. 
One of the main contributions of this paper is the 
classification of the geometric and topological 
configurations arising in the equality case, together 
with a description of the corresponding restrictions 
on the second fundamental form. 

For $n\geq5$, the classification includes spherical 
space forms, manifolds whose universal cover splits 
isometrically as the product of a line 
and a compact manifold diffeomorphic 
to a sphere, and additional configurations 
arising in the equality case. In the latter setting, 
the geometry of the irreducible factors of the universal 
cover is restricted, while the immersion is either 
totally geodesic or exhibits a specific structure 
involving the relative nullity and a principal normal 
of high multiplicity. In dimension four, further 
phenomena arise, including K\"ahler structures, 
Riemannian products of surfaces, and flat manifolds.

The proofs combine a detailed analysis of the equality 
cases in the algebraic estimates for isotropic 
curvature with global results on manifolds with 
nonnegative isotropic curvature and 
Bochner--Weitzenb\"ock techniques. In dimension four, 
the special structure of two-forms plays a crucial 
role.

Before stating our main theorems, we recall that a vector 
$\delta\in N_fM(x)$ is called a principal normal of 
$f$ at $x$ if the subspace 
$$
E_\delta(x)
=
\left\{
X\in T_xM:
\alpha_f(X,Y)=\langle X,Y\rangle\delta
\quad\text{for all}\quad Y\in T_xM
\right\}
$$
is nontrivial (see, for instance, \cite{DT}). 
The dimension of $E_\delta(x)$ is called the 
multiplicity of $\delta$. A principal normal 
is said to be a Dupin principal normal if its 
multiplicity is at least two. 
The relative nullity subspace of $f$ at $x$ is 
defined by 
$$
\mathcal D_f(x)
=
\left\{
X\in T_xM:
\alpha_f(X,Y)=0
\quad\text{for all}\quad Y\in T_xM
\right\}.
$$
Throughout, all submanifolds are assumed 
to be connected. 

The first main theorem establishes the 
higher-dimensional classification.

\begin{theorem}\label{n5}
Let $f\colon M^n\to\bar M^{n+m}$,
$n\geq5$, be an isometric immersion of a compact
Riemannian manifold. Suppose that \eqref{1} holds
at every point, where $c\geq0$ is a lower bound on
the curvature operator of $\bar M^{n+m}$.
Then one of the following alternatives holds:
\begin{enumerate}[
label=(\roman*),
topsep=1pt,
itemsep=0.2ex,
leftmargin=*,
partopsep=0.2ex,
parsep=0.2ex
]

\item The manifold $M^n$ is
diffeomorphic to a spherical space form.

\item The universal cover of $M^n$ is
isometric to a Riemannian product
$\R\times N^{n-1}$, where $N^{n-1}$
is diffeomorphic to $\Sf^{n-1}$ and
has nonnegative isotropic curvature.

\item Equality holds in \eqref{1} at
every point, $c=0$, and one of the
following occurs:
\begin{enumerate}[
label=(\alph*),
leftmargin=*,
itemsep=0.2ex,
partopsep=0.2ex,
parsep=0.2ex
]

\item The immersion $f$ is totally geodesic,
and the universal cover $\tilde M^n$ of $M^n$
is isometric to a Riemannian product
\[
(\R^{n_0},g_0)\times
(N_1^{n_1},g_1)\times\cdots\times
(N_r^{n_r},g_r),
\]
where $n_0\geq0$, $r\geq0$, $g_0$ is
the Euclidean metric, and each
$(N_i,g_i)$, $1\leq i\leq r$,
is a compact, simply connected, irreducible
Riemannian manifold of dimension $n_i\geq2$
with nonnegative curvature operator.
The Euclidean factor is omitted when $n_0=0$,
and the compact factors are omitted when $r=0$.

\item The immersion $f$ is not totally geodesic,
and the universal cover $\tilde M^n$
is isometric to a Riemannian product
$(P^2,g_P)\times(N^{n-2},g_N)$,
where $(P^2,g_P)$ is either compact
or isometric to the Euclidean plane,
and $(N^{n-2},g_N)$ is compact,
simply connected, irreducible,
and has nonnegative curvature operator.
Moreover, at every point where $f$ is not
totally geodesic, the normal vector
$\delta=n\mathcal H_f/(n-2)$ is a Dupin
principal normal of $f$ of multiplicity $n-2$ 
and we have $\pi_*(TP)=\mathcal D_f$ and 
$\pi_*(TN)=E_\delta$, 
where $\pi\colon\tilde M^n\to M^n$
is the universal covering map.
\end{enumerate}
\end{enumerate}

In alternative~(iii), every compact
irreducible factor of dimension at least
three is either a compact, simply connected,
irreducible symmetric space,
a Riemannian manifold diffeomorphic
to a sphere, %and admitting a metric with positive curvature operator
or a K\"ahler
manifold biholomorphic to a complex
projective space. Every compact
two-dimensional factor is diffeomorphic
to $\Sf^2$ and has nonnegative
Gaussian curvature.
\end{theorem}

The four-dimensional case requires a separate analysis, 
owing to the special geometric and topological phenomena 
that arise in this dimension. In particular, the equality 
case permits K\"ahler and product structures for which 
more precise information about the extrinsic geometry 
can be obtained. Recall that an isometric immersion 
of a K\"ahler manifold is called $(2,0)$-geodesic if 
the $(2,0)$-component of its second fundamental form 
vanishes. This notion was introduced by Eschenburg 
and Tribuzy \cite{ET}; see also \cite{ET1,eft}.

\begin{theorem}\label{n=4}
Let $f\colon M^4\to\bar M^{4+m}$ 
be an isometric immersion of a compact, 
oriented Riemannian four-manifold. Suppose that 
inequality \eqref{1} holds at every point, where 
$c\geq0$ is a lower bound on the curvature 
operator of $\bar M^{4+m}$. 
Let $\pi\colon\tilde M^4\to M^4$ denote 
the universal covering map. 
Then one of the following assertions holds: 
\begin{enumerate}[
label=(\roman*),
topsep=1pt,
itemsep=0.2ex,
leftmargin=*,
partopsep=0.2ex,
parsep=0.2ex,
]

\item The manifold $M^4$ is diffeomorphic 
to $\Sf^4$. 

\item The universal cover $\tilde M^4$ is 
isometric to a Riemannian product 
$\R\times N^3$, where $N^3$ is 
diffeomorphic to $\Sf^3$ and has 
nonnegative Ricci curvature. 

\item Equality holds in \eqref{1} 
at every point, and one of the following 
occurs: 
\begin{enumerate}[
label=(\alph*),
topsep=1pt,
itemsep=0.2ex,
partopsep=0.2ex,
parsep=0.2ex,
leftmargin=*
]
\item The manifold $M^4$ admits a K\"ahler 
structure with respect to which $f$ is 
$(2,0)$-geodesic. Moreover, $M^4$ is 
biholomorphic to the complex projective 
plane $\CP^2$. 

\item The universal cover $\tilde M^4$ 
is isometric to a Riemannian product 
$(N_1,g_1)\times(N_2,g_2)$, 
where each factor is diffeomorphic to 
$\Sf^2$ and has nonnegative Gaussian 
curvature. Moreover, there exist smooth 
Dupin principal normals $\delta_1$ 
and $\delta_2$ of $f$ satisfying 
\[
\pi_*(TN_i)\subset E_{\delta_i},
\quad i=1,2,
\qquad
\langle\delta_1,\delta_2\rangle=-c.
\]

\item We have $c=0$, and 
the universal cover $\tilde M^4$ 
is isometric to a Riemannian product 
$(\R^2,g_0)\times(N^2,g_N)$, 
where $g_0$ is the Euclidean metric 
and $N^2$ is diffeomorphic to $\Sf^2$ 
with nonnegative Gaussian curvature. 
Moreover, $\pi_*(T\R^2)\subset\mathcal D_f$, 
and the normal vector field 
$\delta=2\mathcal H_f$ is a Dupin 
principal normal satisfying 
$\pi_*(TN^2)\subset E_\delta$. 
\item We have $c=0$, the manifold 
$M^4$ is flat with 
$1\leq b_1(M)\leq4$, and 
the immersion $f$ is totally geodesic. 
\end{enumerate} 
\end{enumerate} 
\end{theorem} 

Theorems~\ref{n5} and~\ref{n=4} provide a unified 
description of the intrinsic geometry, topology, 
and extrinsic structure of submanifolds satisfying 
\eqref{1}, well beyond the space-form setting. 
A noteworthy feature of our approach is that the 
ambient curvature lower bound is required only on 
tangential bivectors along the immersion. 
The conclusions of these theorems cannot, however, 
in general be strengthened to extrinsic rigidity 
statements. Indeed, in Section~6 we construct 
families of geometrically distinct immersions 
satisfying \eqref{1}, including immersions of 
manifolds diffeomorphic to $\Sf^n$ and of products 
$\Sf^1\times N^{n-1}$, where $N^{n-1}$ is 
diffeomorphic to $\Sf^{n-1}$ and has positive 
sectional curvature.

We next investigate a stronger pinching condition, 
corresponding to the case $k=1$ of the aforementioned 
family of inequalities. In contrast to \eqref{1}, 
this condition ensures nonnegative sectional 
curvature and yields a more restrictive classification, 
together with sharper conclusions concerning the 
extrinsic geometry of the immersion.

\begin{theorem}\label{k=1}
Let $f\colon M^n\to\bar M^{n+m}$, $n\geq4$,
be an isometric immersion of a compact
Riemannian manifold, and let $c\geq0$
be a lower bound on the curvature operator
of $\bar M^{n+m}$. Suppose that
\begin{equation}\label{2}
S\leq 2c+\frac{n^2H^2}{n-1}
\tag{$**$}
\end{equation}
holds at every point. Then $M^n$ has
nonnegative sectional curvature, and 
one of the following alternatives holds:
\begin{enumerate}[
label=(\roman*),
topsep=1pt,
itemsep=0.5ex,
leftmargin=*,
partopsep=0.2ex,
parsep=0.2ex
]

\item The manifold $M^n$ is diffeomorphic
to a spherical space form (to $\Sf^4$
if $n=4$ and $M^4$ is orientable). 

\item Equality holds in \eqref{2} at
every point, $c=0$, and one of the
following occurs:
\begin{enumerate}[
label=(\alph*),
leftmargin=*,
itemsep=0.2ex,
partopsep=0.2ex,
parsep=0.2ex
]

\item The immersion $f$ is totally geodesic, 
and its universal cover is isometric to 
a Riemannian product 
\[
(\R^{n_0},g_0)\times
(N_1^{n_1},g_1)\times\cdots\times
(N_r^{n_r},g_r),
\]
where $n_0\geq0$, $r\geq0$, $g_0$ is
the Euclidean metric, and each
$(N_i,g_i)$, $1\leq i\leq r$,
is a compact, simply connected, irreducible
Riemannian manifold of dimension $n_i\geq2$
with nonnegative curvature operator.
The Euclidean factor is omitted when $n_0=0$,
and the compact factors are omitted when $r=0$.
Every compact irreducible factor of dimension
at least three is either a compact, simply
connected, irreducible symmetric space,
a Riemannian manifold diffeomorphic to a sphere, 
%and admitting a metric with positive curvature operator, 
or a K\"ahler manifold biholomorphic
to a complex projective space. Every
two-dimensional factor is diffeomorphic to
$\Sf^2$ and has nonnegative Gaussian curvature.

\item The immersion $f$ is not totally
geodesic, and the universal cover of $M^n$
is isometric to a Riemannian product
$\R\times N^{n-1}$, where $N^{n-1}$
is diffeomorphic to $\Sf^{n-1}$ and
has nonnegative curvature operator,
which is positive at some point.
Moreover, at every point where $f$ is not
totally geodesic, the normal vector
$\delta=n\mathcal H_f/(n-1)$ is a Dupin
principal normal of $f$ of multiplicity $n-1$,
and we have $\pi_*(T\R)=\mathcal D_f$ and
$\pi_*(TN)=E_\delta$, where
$\pi\colon\tilde M^n\to M^n$ is the
universal covering map.
\end{enumerate}
\end{enumerate}
\end{theorem}

Together, Theorems~\ref{n5}, \ref{n=4}, 
and~\ref{k=1} provide a detailed description of 
the intrinsic and extrinsic consequences of the 
nonstrict pinching inequalities
\[
S\leq 2kc+\frac{n^2H^2}{n-k},
\qquad 1\leq k\leq n/2,
\]
for $k=1$ and $k=2$. These results suggest the 
more general problem of determining the geometric 
and topological restrictions imposed by the above 
inequalities for the remaining admissible integer 
values of $k$. Of particular interest is the 
classification of the configurations arising in 
the equality case in general ambient manifolds, 
beyond the strict space-form setting considered 
by Howard and Wei \cite{HW}.

In the space-form setting, the higher-dimensional 
classification simplifies considerably: even when 
equality is allowed in the pinching condition, 
only two alternatives remain. The following result, 
previously established in \cite{v5}, is an immediate 
consequence of Theorem~\ref{n5}. 
We denote by $\Q_c^d$ the complete, 
simply connected Riemannian manifold 
of dimension $d$ and constant sectional 
curvature $c$. 

\begin{corollary}\label{sfhigh} 
Let $f\colon M^n\to\Q_c^{n+m}$, with $n\geq5$, 
be an isometric immersion of a compact Riemannian 
manifold into a space form of curvature $c\geq0$. 
Suppose that 
$$
S\leq4c+\frac{n^2H^2}{n-2} 
$$
at every point. Then the manifold $M^n$ is 
diffeomorphic to a spherical space form, 
or its universal cover is isometric to a 
Riemannian product $\R\times N$, where $N$ is 
diffeomorphic to $\Sf^{n-1}$ and has nonnegative 
isotropic curvature. 
\end{corollary}

The four-dimensional classification in space forms, 
previously established in Theorem~2 of \cite{v4}, 
can also be recovered by combining 
Theorem~\ref{n=4} with Theorem~1.1 of \cite{eft} 
and Proposition~12 of \cite{v3}.

\begin{corollary}\label{sf4}
Let $f\colon M^4\to\Q_c^{4+m}$, $c\geq0$, be 
an isometric immersion of a compact, oriented 
Riemannian four-manifold. Suppose that \eqref{1} 
is satisfied at every point. Then one of the following 
holds: 
\begin{enumerate}[
label=(\roman*),
topsep=1pt,
itemsep=0.2ex,
leftmargin=*,
partopsep=0.2ex,
parsep=0.2ex,
]

\item $M^4$ is diffeomorphic to $\Sf^4$. 
\item The universal cover of $M^4$ is 
isometric to a Riemannian product $\R\times N$, 
where $N$ is diffeomorphic to $\Sf^3$ 
and has nonnegative Ricci curvature. 
\item Equality holds in \eqref{1} at every point, 
and one of the following occurs: 
\begin{enumerate}[
label=(\alph*),
topsep=1pt,
itemsep=0.2ex,
leftmargin=*,
partopsep=0.2ex,
parsep=0.2ex,
]

\item $M^4$ is isometric to a Riemannian product 
$\Sf^2(r)\times\Sf^2\big(\sqrt{R^2-r^2}\big)$, 
and $f$ is a composition $f=j\circ g$, where 
$g\colon M^4\to \Sf^5(R)$ 
is the standard embedding of 
$\Sf^2(r)\times\Sf^2(\sqrt{R^2-r^2})$ into $\Sf^5(R)$ 
and $j\colon\Sf^5(R)\to\Q_c^{4+m}$ 
is an umbilical inclusion. 

\item $M^4$ is isometric to the complex projective plane 
$\CP_r^2$ of constant holomorphic sectional curvature $4/(3r^2)$, 
where $r=1/\sqrt{c+H^2}$, and $f=j\circ g$. Here 
$g\colon\CP_r^2\to\Sf^7(r)$ is the standard 
minimal embedding, and 
$j\colon\Sf^7(r)\to\Q_c^{4+m}$ 
is an umbilical inclusion. 
\end{enumerate} 
\end{enumerate} 
\end{corollary} 

As a further consequence of Theorem~\ref{k=1}, 
the stronger pinching condition yields a 
differentiable sphere theorem for submanifolds 
of space forms, even when equality is allowed. 
This complements the corresponding strict 
pinching result of Howard and Wei \cite{HW}; 
see also \cite{XG}.

\begin{corollary}\label{cor}
Let $f\colon M^n\to\Q_c^{n+m}$, with $n\geq4$ and $c\geq0$, 
be an isometric immersion of a compact, oriented Riemannian 
manifold. If 
$$
S\leq2c+\frac{n^2H^2}{n-1} 
$$
at every point, then $M^n$ is diffeomorphic to $\Sf^n$. 
\end{corollary} 

The paper is organized as follows. Section~2 develops 
the preliminaries on isotropic curvature and establishes 
the local algebraic results underlying the analysis of 
the equality cases. It also contains results on 
immersions of Riemannian products, Bochner--Weitzenb\"ock 
techniques, and applications in even dimensions. 
Section~3 establishes the four-dimensional 
classification in Theorem~\ref{n=4}, while Section~4 
proves Theorem~\ref{n5}, treating the cases of finite 
and infinite fundamental groups separately. Section~5 
is devoted to the proof of Theorem~\ref{k=1} and its 
applications to submanifolds of space forms. Finally, 
Section~6 presents examples illustrating the scope 
and sharpness of our results.

\section{Preliminaries and auxiliary results} 
\subsection{Isotropic curvature} 
In this subsection, we recall the notions of positive and 
nonnegative isotropic curvature introduced by 
Micallef and Moore in \cite{MM}. 

Let $(M^n,g)$ be a Riemannian manifold of dimension 
$n\geq4$ with Levi-Civit\'a connection 
$\n$ and curvature tensor $R$ given by 
$$
R(X,Y)=[\n_X,\n_Y]-\n_{[X,Y]},\qquad X,Y\in\mathcal X(M). 
$$
Recall that the \emph{curvature operator} at a point 
$x\in M^n$ is the self-adjoint endomorphism 
$$
\mathcal R\colon\Lambda^2T_xM\to\Lambda^2T_xM 
$$
on the second exterior power $\Lambda^2T_xM$ of 
the tangent space, defined by 
$$
\<\<\mathcal R(X\wedge Y),Z\wedge W\>\>
=\langle R(X,Y)W,Z\rangle,\qquad X,Y,Z,W\in T_xM, 
$$
where the inner products $\<\<\cdot,\cdot\>\>$ and 
$\<\cdot,\cdot\>$ are induced by the metric $g$. 

The inner product $\<\cdot,\cdot\>$ on $T_xM$ can 
be extended to the complexification $T_xM\otimes\C$ 
either as a complex bilinear form $(\cdot,\cdot)$ or as a 
Hermitian inner product $((\cdot,\cdot))$ defined by 
$$
((v,w))=(v,\overline w),\qquad v,w\in T_xM\otimes\C.
$$
The curvature operator is also extended to a 
complex-linear map
$$
\mathcal R\colon\Lambda^2T_xM\otimes\C
\to\Lambda^2T_xM\otimes\C.
$$
For a complex two-plane $\sigma\subset T_xM\otimes\C$, 
its \emph{complex sectional curvature} 
$K^c(\sigma)$ is defined by 
$$
K^c(\sigma)=
\frac{((\mathcal R(v\wedge w),v\wedge w))}{\|v\wedge w\|^2}, 
$$
where $\{v,w\}$ is a basis of $\sigma$ and the norm 
is induced by the Hermitian inner product. A two-plane 
$\sigma\subset T_xM\otimes\C$ is called 
\emph{isotropic} if $(v,v)=0$ for all $v\in\sigma$. 

We say that $(M^n,g)$ has \emph{positive isotropic curvature} 
(respectively, \emph{nonnegative isotropic curvature}) 
at a point $x\in M^n$ if $K^c(\sigma)>0$ (respectively, 
$K^c(\sigma)\geq0$) for every isotropic two-plane at 
$x$. 

To each orthonormal four-frame 
$F=\{e_1,e_2,e_3,e_4\}\subset T_xM$ at a point $x\in M^n$, 
we associate the isotropic two-plane 
$$
\sigma_F=\spa_\C\left\{e_1+\mathsf ie_2,e_3+\mathsf ie_4\right\}, 
$$
where $\mathsf i$ is the imaginary unit. It follows 
(see \cite{MM}) that 
$$
K^c(\sigma_F)=\frac{1}{4}\bigl(R_{1313} + R_{1414}
+ R_{2323} + R_{2424}-2 R_{1234}\bigr),
$$
where the components of the curvature 
tensor $R$ are given by 
$$
R_{ijkl}=g\big(R(e_i,e_j)e_l, e_k\big), 
\qquad 1\leq i,j,k,l\leq4. 
$$ 
Hence, $(M^n,g)$ has positive 
isotropic curvature (respectively, nonnegative isotropic 
curvature) at a point $x\in M^n$ if 
$$
R_{1313} + R_{1414}+ R_{2323} + R_{2424}
-2 R_{1234}>0, 
$$
(respectively, 
$$
R_{1313} + R_{1414}+ R_{2323} + R_{2424}
-2 R_{1234} \geq0), 
$$
for all orthonormal four-frames 
$\{e_i\}_{i=1}^4\subset T_xM$. 
The manifold $(M^n,g)$ is said to have nonnegative 
(or positive) isotropic curvature if it satisfies the 
corresponding condition at every point and for all 
orthonormal four-frames. We follow 
this notation throughout the paper. 

It is well known that symmetric spaces of compact type 
have nonnegative isotropic curvature. Moreover, if such 
a space has nonconstant sectional curvature, then at 
each point there exists an orthonormal four-frame 
with vanishing complex sectional curvature. 

Let $M$ be a symmetric space with $\dim M\geq4$. 
For each $x\in M$, let $\mathcal F_x$ denote the set 
of all orthonormal four-frames 
$F=\{e_i\}_{i=1}^4\subset T_xM$ with vanishing 
complex sectional curvature. Each $F\in\mathcal F_x$ 
spans a four-dimensional subspace 
$V_F=\spa\left\{e_i\right\}_{i=1}^4\subset T_xM$. 
The following result, whose proof is given in 
the appendix of \cite{v5}, will be used later. 

\begin{proposition}\label{fillup}
Let $M$ be a compact, simply connected, irreducible 
symmetric space of dimension at least five with 
nonconstant sectional curvature. Then, at every 
point $x\in M$, the orthogonal complements of the 
subspaces $V_F$, $F\in\mathcal F_x$, span $T_xM$; 
that is,
$$
\sum_{F\in\mathcal F_x}V_F^\perp=T_xM.
$$
\end{proposition}

We next recall two classification results for 
compact Riemannian manifolds with nonnegative 
isotropic curvature that will be used in the sequel. 
Recall that a Riemannian manifold is said to be 
\emph{locally irreducible} if it does not admit 
a nontrivial local Riemannian product decomposition. 

\begin{theorem}[\cite{MW}]\label{duke}
Let $(M^n,g)$, $n\geq4$, be a compact, locally 
reducible Riemannian manifold with nonnegative 
isotropic curvature. Then its Riemannian universal 
cover $(\tilde M^n,\tilde g)$ is isometric to one 
of the following Riemannian products: 
\begin{enumerate}[
label=(\alph*),
topsep=1pt,
itemsep=0.2ex,
leftmargin=*,
partopsep=0.2ex,
parsep=0.2ex
]
\item
$(\R^{n_0},g_0)\times
(N_1^{n_1},g_1)\times\cdots\times
(N_r^{n_r},g_r)$,
where $n_0\geq0$, $g_0$ is the Euclidean metric, 
and, for each $1\leq i\leq r$, either 
$n_i=2$ and $N_i$ is diffeomorphic to $\Sf^2$ 
with nonnegative Gaussian curvature, or 
$n_i\geq3$ and $(N_i,g_i)$ is a compact, 
irreducible Riemannian manifold with 
nonnegative flag curvature. 

\item
$(\varSigma^2,g_\varSigma)\times
(N^{n-2},g_N)$,
where $\varSigma^2$ is a surface with Gaussian 
curvature $K_\varSigma$, and $(N,g_N)$ is a 
compact, irreducible Riemannian manifold whose 
sectional curvature $K_N$ satisfies 
$\min K_N\geq-\min K_{\varSigma}>0$. 
\end{enumerate}
\end{theorem}

\begin{theorem}[\cite{BSacta}]\label{acta}
Let $(M^n,g)$, $n\geq4$, be a compact, locally 
irreducible Riemannian manifold with nonnegative 
isotropic curvature. Then one of the following 
alternatives holds: 
\begin{enumerate}[
label=(\roman*),
topsep=1pt,
itemsep=0.2ex,
leftmargin=*,
partopsep=0.2ex,
parsep=0.2ex
]
\item
$M^n$ is diffeomorphic to a spherical space form. 

\item
The universal cover of $M^n$ is a K\"ahler 
manifold biholomorphic to $\CP^{n/2}$. 

\item
The universal cover of $M^n$ is isometric 
to a compact symmetric space. 
\end{enumerate}
\end{theorem}

\subsection{The pinching condition and isotropic curvature}
The aim of this subsection is to investigate in detail the relationship 
between the pinching condition \eqref{1} and isotropic curvature. 
The auxiliary results obtained are pointwise in nature. 

The following algebraic lemma was proved in \cite[Lemma 10]{v5}. 

\begin{lemma}\label{alg2}
Let $f\colon M^n\to\bar {M}^{n+m}$, 
$n\geq4$, be an isometric immersion. 
For every point $x\in M^n$, every orthonormal basis 
$\{e_i\}_{i=1}^n$ of $T_xM$, and every 
orthonormal basis $\{\xi_\alpha\}_{\alpha=1}^m$ 
of $N_fM(x)$, the following statements hold: 
\begin{enumerate}[
label=(\roman*),
topsep=1pt,
itemsep=0.2ex,
leftmargin=*,
partopsep=0.2ex,
parsep=0.2ex,
]

\item For each $\varepsilon=\pm1$, we have 
\begin{align}\label{ineq2alg}
\sum_{\alpha=1}^m\Big(&\sum_{i=1}^{2}\sum_{j=3}^{4}\big(h^{\alpha}_{ii}h^{\alpha}_{jj}-(h^\alpha_{ij})^2\big)+2\varepsilon(h^{\alpha}_{14}h^{\alpha}_{23}-h^{\alpha}_{13}h^{\alpha}_{24})\Big)
\geq \frac{n^2H^2}{n-2}-S, 
\end{align}
where $h^\alpha_{ij}=\langle A_{\xi_\alpha}e_i,e_j\rangle$ 
and $A_{\xi_\alpha}$ denotes the shape operator 
associated with $\xi_\alpha$. 
\item Suppose that equality holds in 
\eqref{ineq2alg} for some $\varepsilon=\pm1$. 
Then, for every $1\leq\alpha\leq m$, 
\begin{equation}\label{ineq3alg}
\begin{gathered}
h^\alpha_{11}=h^\alpha_{22},\quad
h^\alpha_{33}=h^\alpha_{44},\quad
h^\alpha_{12}=h^\alpha_{34}=0,
\quad h^\alpha_{23}=-\varepsilon h^\alpha_{14},
\quad
h^\alpha_{24}=\varepsilon h^\alpha_{13}.
\end{gathered}
\end{equation}
Moreover, if $n\geq5$, there exist real 
numbers $\rho_\alpha$, $1\leq\alpha\leq m$, 
such that 
\begin{equation}\label{ineq4alg}
h^\alpha_{ij}=0\quad\text{for all }1\leq i\leq4,\; 5\leq j\leq n,
\quad h^\alpha_{ij}=\delta_{ij}\rho_\alpha
\quad\text{for all }5\leq i,j\leq n,
\end{equation}
and 
\begin{equation}\label{ineq5alg}
h^\alpha_{11}+h^\alpha_{44}=\rho_\alpha,
\qquad 1\leq\alpha\leq m.
\end{equation}
\end{enumerate}
\end{lemma}
\vspace{1ex}

Let $f\colon M^n\to\bar {M}^{n+m}$, 
$n\geq4$, be an isometric immersion of a 
Riemannian manifold $(M^n,\langle\cdot,\cdot\rangle)$ 
into a Riemannian manifold 
$(\bar {M}^{n+m},\langle\cdot,\cdot\rangle)$. 
Denote by $R_M$ and $R_{\bar M}$ 
the curvature tensors of $M$ and 
$\bar {M}$, respectively. 

The \emph{restricted curvature tensor} of $f$ 
is the tensor field 
\[
R_f\colon
TM\times TM\times TM\to TM
\]
defined by 
\[
\langle R_f(X,Y)Z,W\rangle
=
\langle R_{\bar M}(f_*X,f_*Y)f_*Z,
f_*W\rangle,
\]
for all $X,Y,Z,W\in TM$. 
The Gauss equation for $f$ takes the form 
\begin{align*}
\langle R_M(X,Y)Z,W\rangle
={}&
\langle R_f(X,Y)Z,W\rangle
+\langle\alpha_f(X,W),\alpha_f(Y,Z)\rangle
-\langle\alpha_f(X,Z),\alpha_f(Y,W)\rangle,
\end{align*}
where $\alpha_f$ is the second fundamental 
form of $f$. 

The \emph{restricted curvature operator} of $f$ 
is the self-adjoint bundle endomorphism 
\[
{\mathcal R}_f\colon
\Lambda^2TM\to\Lambda^2TM
\]
defined by 
\[
\langle\!\langle
{\mathcal R}_f(X\wedge Y),
Z\wedge W
\rangle\!\rangle
=
\langle
R_f(X,Y)W,Z
\rangle,
\]
for all $X,Y,Z,W\in TM$. 
The Gauss equation then implies that 
the curvature operator of $M^n$ is 
given by 
\begin{align}
\langle\!\langle
\mathcal R_M(X\wedge Y),Z\wedge W
\rangle\!\rangle
={}&
\langle\!\langle
{\mathcal R}_f(X\wedge Y),
Z\wedge W
\rangle\!\rangle\nonumber
+\langle\alpha_f(X,Z),\alpha_f(Y,W)\rangle
\nonumber\\
&-\langle\alpha_f(X,W),\alpha_f(Y,Z)\rangle.
\label{RR}
\end{align}

This notation will be used throughout the paper. 

\begin{lemma}\label{algc}
Let $f\colon M^n\to\bar {M}^{n+m}$, $n\geq4$, 
be an isometric immersion into a Riemannian 
manifold $\bar {M}^{n+m}$ whose curvature 
operator is bounded below by a constant $c$. 
Then, at every point $x\in M^n$ and for every 
orthonormal four-frame $\{e_i\}_{i=1}^4\subset T_xM$, 
the following assertions hold at $x$: 
\begin{enumerate}[
label=(\roman*),
topsep=1pt,
itemsep=0.2ex,
leftmargin=*,
partopsep=0.2ex,
parsep=0.2ex,
]

\item For each $\varepsilon=\pm1$, we have 
\begin{align}\label{ineq1}
R_{1313}+R_{1414}+R_{2323}+R_{2424}-2\varepsilon R_{1234}
\geq4c+\frac{n^2H^2(x)}{n-2}-S(x). 
\end{align}

\item Suppose that equality holds in 
\eqref{ineq1} for some $\varepsilon=\pm1$ 
and an orthonormal four-frame 
$F=\{e_i\}_{i=1}^4\subset T_xM$. 
Then 
\begin{align}
\alpha_{11}=\alpha_{22},\quad\alpha_{33}
=\alpha_{44},\quad\alpha_{12}=\alpha_{34}=0,\quad 
\alpha_{23}=-\varepsilon &\alpha_{14},\quad\alpha_{24}
=\varepsilon\alpha_{13}, \label{c1}\\ 
\alpha_{11}+\alpha_{44}=\frac{n}{n-2}\mathcal H_f(x), \label{c2}
\end{align} 
where, for convenience, we write 
$\alpha_{ij}=\alpha_f(e_i,e_j)$. 
Moreover, if $n\geq5$, the normal vector 
$\delta=n\mathcal H_f/(n-2)$ is a 
principal normal of $f$ satisfying 
$V_F^\perp\subset E_\delta(x)$, where 
$V_F^\perp$ denotes the orthogonal complement 
of $V_F=\spa\left\{e_i\right\}_{i=1}^4$ in 
$T_xM$. 
\end{enumerate}
\end{lemma}
\proof
By direct computation, we obtain 
\begin{align*}
\langle\!\langle \mathcal R_MV,V\rangle\!\rangle
+\langle\!\langle \mathcal R_MW,W\rangle\!\rangle
=R_{1313}+R_{1414}+R_{2323}+R_{2424}
-2\varepsilon(R_{1324}-R_{1423}),
\end{align*}
where 
$$ 
V=e_1\wedge e_3-\varepsilon e_2\wedge e_4,\qquad
W=e_1\wedge e_4+\varepsilon e_2\wedge e_3.
$$
Using the algebraic symmetries of the curvature tensor 
and the first Bianchi identity, we conclude that 
\begin{equation}\label{Rf1}
\begin{aligned}
\langle\!\langle \mathcal R_MV,V\rangle\!\rangle
+\langle\!\langle \mathcal R_MW,W\rangle\!\rangle
=R_{1313}+R_{1414}+R_{2323}+R_{2424}
-2\varepsilon R_{1234}.
\end{aligned}
\end{equation}

Setting 
\[
\bar V=f_*e_1\wedge f_*e_3
-\varepsilon f_*e_2\wedge f_*e_4,
\qquad
\bar W=f_*e_1\wedge f_*e_4
+\varepsilon f_*e_2\wedge f_*e_3,
\]
we have 
\be\label{Rf2}
\<\<{\mathcal R}_fV,V\>\>+\<\<{\mathcal R}_fW,W\>\>=
\<\<{\mathcal R}_{\bar M}\bar V,\bar V\>\>+
\<\<{\mathcal R}_{\bar M}\bar W,\bar W\>\>
\ge4c. 
\ee

Extend $F$ to an orthonormal basis 
$\{e_i\}_{i=1}^n$ of $T_xM$ and choose 
an orthonormal basis 
$\{\xi_\alpha\}_{\alpha=1}^m$ of $N_fM(x)$. 
Using \eqref{RR}, \eqref{Rf1}, \eqref{Rf2}, 
and part $(i)$ of Lemma \ref{alg2}, we obtain 
\begin{align}\label{algc1}
R_{1313}&+R_{1414}+R_{2323}+R_{2424}-2\varepsilon R_{1234}\nonumber\\
&=\<\<\mathcal R_MV,V\>\>+\<\<\mathcal R_MW,W\>\>
=\<\<\mathcal R_fV,V\>\>+\<\<\mathcal R_fW,W\>\>\nonumber\\
&+\sum_{\alpha=1}^m\Big(\sum_{i=1}^{2}\sum_{j=3}^{4}\big(h^{\alpha}_{ii}h^{\alpha}_{jj}-(h^\alpha_{ij})^2\big)+2\varepsilon(h^{\alpha}_{14}h^{\alpha}_{23}-h^{\alpha}_{13}h^{\alpha}_{24})\Big)\nonumber\\
&\geq4c+\frac{n^2H^2(x)}{n-2}-S(x). 
\end{align}
This proves \eqref{ineq1}. 

Assume now that equality holds in \eqref{ineq1} 
for an orthonormal four-frame 
$F=\{e_i\}_{i=1}^4\subset T_xM$. 
Then, all inequalities used in \eqref{algc1} 
must be equalities. In particular, equality 
holds in \eqref{ineq2alg} and \eqref{Rf2}. 
Hence, part $(ii)$ of Lemma \ref{alg2}
applies, and \eqref{c1} follows from
\eqref{ineq3alg}.

If $n=4$, the first two equalities in 
\eqref{c1} yield 
$\alpha_{11}+\alpha_{44}
=2\mathcal H_f(x)$,
which is precisely \eqref{c2}. 

If $n\geq5$, it follows from 
\eqref{ineq4alg} that the normal vector 
$\delta=\sum_{\alpha=1}^m\rho_\alpha\xi_\alpha$ 
is a principal normal of $f$ satisfying 
$V_F^\perp\subset E_\delta(x)$. 
Moreover, \eqref{ineq5alg} gives 
$\alpha_{11}+\alpha_{44}=\delta$. 
Taking the trace of the second fundamental 
form, we obtain 
$n\mathcal H_f(x)=(n-2)\delta$. 
This establishes \eqref{c2} and completes 
the proof. 
\qed
\vspace{1ex}

The following proposition, which establishes 
the relationship between the inequality \eqref{1} 
and the isotropic curvature, is one of the key 
auxiliary results used in our proofs. 

\begin{proposition}\label{propu}
Let $f\colon M^n\to\bar {M}^{n+m}$, $n\geq4$, 
be an isometric immersion into a Riemannian 
manifold $\bar {M}^{n+m}$ whose curvature operator 
is bounded below by a constant $c$. Suppose that the 
inequality \eqref{1} holds at a point $x\in M^n$. Then the 
following statements hold at $x$: 
\begin{enumerate}[
label=(\roman*),
topsep=1pt,
itemsep=0.2ex,
leftmargin=*,
partopsep=0.2ex,
parsep=0.2ex,
]

\item The manifold $M^n$ has nonnegative 
isotropic curvature at $x$. 
\item Suppose that $M^n$ does not have 
positive isotropic curvature at $x$. Then equality 
holds in \eqref{1} at $x$. Moreover, for any orthonormal 
four-frame $F=\{e_i\}_{i=1}^4\subset T_xM$ satisfying 
\be\label{eq1}
R_{1313}+R_{1414}+R_{2323}+R_{2424}
-2\varepsilon R_{1234}=0,
\ee
with $\varepsilon=\pm1$, the conditions 
\eqref{c1}--\eqref{c2} in part $(ii)$ of 
Lemma \ref{algc} are satisfied at $x$. 
If $n\geq5$, the normal vector 
$\delta=n\mathcal H_f/(n-2)$ 
is a principal normal of $f$ satisfying 
$V_F^\perp\subset E_\delta(x)$, where 
$V_F=\spa\left\{e_i\right\}_{i=1}^4$. 
\end{enumerate}
\end{proposition}
\proof 
\noindent $(i)$ 
That $M^n$ has nonnegative isotropic 
curvature at $x$ follows directly from 
\eqref{1} and part $(i)$ of Lemma \ref{algc}. 
\vspace{1ex}

\noindent $(ii)$ 
Suppose that $M^n$ does not have 
positive isotropic curvature at $x$. 
Then there exists an orthonormal four-frame 
$F=\{e_i\}_{i=1}^4\subset T_xM$ satisfying 
\eqref{eq1} for some $\varepsilon=\pm1$. 
By \eqref{1} and \eqref{ineq1}, we have 
\[
0\geq 4c+\frac{n^2H^2(x)}{n-2}-S(x)
\geq 0.
\]
Thus equality holds in \eqref{1} at $x$, 
and part $(ii)$ of Lemma \ref{algc} 
implies that the conditions 
\eqref{c1}--\eqref{c2} are satisfied. 
Moreover, if $n\geq5$, the normal vector 
$\delta=n\mathcal H_f(x)/(n-2)$ is a 
principal normal of $f$ satisfying 
$V_F^\perp\subset E_\delta(x)$. 
\qed
%\vspace{-5pt}
\subsection{Isometric immersions of Riemannian products}

Let $f\colon M^n\to\bar {M}^{n+m}$, 
$n\geq4$, be an isometric immersion of a 
Riemannian manifold $M^n$ isometric to a 
Riemannian product 
\[
(N_1^{n_1},g_1)\times\cdots\times
(N_r^{n_r},g_r),\qquad r\geq2.
\]
Fix a point 
$\hat x=(\hat x_1,\dots,\hat x_r)\in M^n$. 
For each $1\leq i\leq r$, let 
$\sigma_i\colon N_i\to M^n$ denote 
the totally geodesic inclusion defined by 
\[
\sigma_i(x_i)=
(\hat x_1,\dots,x_i,\dots,\hat x_r),
\qquad x_i\in N_i.
\]

Consider the isometric immersions 
\[
f_i=f\circ\sigma_i\colon
N_i\to\bar {M}^{n+m},
\qquad 1\leq i\leq r.
\]
Since each $\sigma_i$ is totally geodesic, 
we have 
\begin{equation}\label{ai}
\alpha_{f_i}(X_i,Y_i)
=
\alpha_f(\sigma_{i*}X_i,\sigma_{i*}Y_i),
\qquad X_i,Y_i\in TN_i.
\end{equation}
Let $H_i$ and 
$\Phi_i=\alpha_{f_i}-g_i\mathcal H_{f_i}$
denote the mean curvature and the 
traceless second fundamental form of 
$f_i$, respectively. 

We shall need the following auxiliary result. 

\begin{lemma}\label{lemmaprod}
Let $f\colon M^n\to\bar {M}^{n+m}$, $n\geq4$, 
be an isometric immersion of a Riemannian product as 
above. Suppose that the curvature operator of 
$\bar {M}^{n+m}$ is bounded below 
by a constant $c$, and that $f$ satisfies the inequality 
\eqref{1} at a point $x=(x_1,\dots,x_r)\in M^n$. 
For each $1\leq i\leq r$, let $f_i=f\circ\sigma_i$ 
be the immersion defined above, with $\hat x=x$.
Then the following inequality holds at $x$: 
\begin{align}\label{pro}
\sum_{i=1}^r\Big(\|\Phi_i\|^2+n_iH_i^2\big(1-\frac{n_i}{n-2}\big)\Big)&+\frac{2(n-3)}{n-2}\sum_{1\leq i<j\leq r}\sum_{k=1}^{n_i}
\sum_{\ell=1}^{n_j}\Big\|\a_f(\sigma_i{_*}e^i_k,\sigma_j{_*}e^j_\ell)\Big\|^2\nonumber\\
&\leq
\frac{c}{n-2}\Big(\sum_{i=1}^rn_i^2-n^2
+4n-8\Big),
\end{align} 
where $\{e^i_k\}_{k=1}^{n_i}$ is an orthonormal basis 
of $T_{x_i}N_i$, $1\leq i\leq r$. 
\end{lemma}
\proof
Fix $x=(x_1,\dots,x_r)\in M^n$ and consider 
the totally geodesic inclusions $\sigma_i$ 
defined using $\hat x=x$. 
Let $\{e_k^i\}_{k=1}^{n_i}$ be an orthonormal 
basis of $T_{x_i}N_i$, $1\leq i\leq r$. 
Then 
\[
\left\{\sigma_{i*}e_k^i:
1\leq i\leq r,\ 1\leq k\leq n_i\right\}
\]
is an orthonormal basis of $T_xM$. 

Using \eqref{ai}, we obtain 
\begin{align}\label{S}
S={}&\sum_{i=1}^r S_i
+2\sum_{1\leq i<j\leq r}
\sum_{k=1}^{n_i}\sum_{\ell=1}^{n_j}
\Big\|
\alpha_f(\sigma_{i*}e_k^i,
\sigma_{j*}e_\ell^j)
\Big\|^2,
\end{align}
where $S_i$ denotes the squared norm of 
the second fundamental form of $f_i$. 
Moreover, 
\[
n\mathcal H_f=\sum_{i=1}^r
n_i\mathcal H_{f_i},
\]
and hence 
\begin{align}\label{H}
n^2H^2={}&\sum_{i=1}^r n_i^2H_i^2
+2\sum_{1\leq i<j\leq r}
n_in_j
\langle\mathcal H_{f_i},
\mathcal H_{f_j}\rangle.
\end{align}

Since $\|\Phi_i\|^2=S_i-n_iH_i^2$, 
it follows from \eqref{S} and \eqref{H} 
that inequality \eqref{1} is equivalent to 
\begin{align}\label{*}
\sum_{i=1}^r\Bigl(\|\Phi_i\|^2+n_iH_i^2\bigl(1-\frac{n_i}{n-2}\bigr)\Bigr)
&+2\sum_{1\leq i<j\leq r}\sum_{k=1}^{n_i}
\sum_{\ell=1}^{n_j}\Bigl\|\a_f(\sigma_i{_*}e^i_k,\sigma_j{_*}e^j_\ell)\Bigr\|^2\nonumber\\
&\leq4c+
\frac{2}{n-2}\sum_{1\leq i<j\leq r}n_in_j\<\mathcal H_{f_i}, \mathcal H_{f_j}\>.
\end{align} 

Since $M$ is a Riemannian product, 
the sectional curvature of every mixed 
two-plane vanishes. Thus, the Gauss 
equation gives 
$$
0=\bigl\langle
{R}_f
(\sigma_{i*}e_k^i,\sigma_{j*}e_\ell^j)
\sigma_{j*}e_\ell^j,
\sigma_{i*}e_k^i
\bigr\rangle+
\langle
\alpha_{f_i}(e_k^i,e_k^i),
\alpha_{f_j}(e_\ell^j,e_\ell^j)
\rangle-
\Big\|
\alpha_f(\sigma_{i*}e_k^i,
\sigma_{j*}e_\ell^j)
\Big\|^2$$
for $i\neq j$. Summing over $k$ and $\ell$, 
and using the lower bound on the ambient 
curvature operator, we obtain 
$$
n_in_j
\langle\mathcal H_{f_i},
\mathcal H_{f_j}\rangle
\leq
-cn_in_j+
\sum_{k=1}^{n_i}\sum_{\ell=1}^{n_j}
\Big\|
\alpha_f(\sigma_{i*}e_k^i,
\sigma_{j*}e_\ell^j)
\Big\|^2.
$$

Substituting this inequality into \eqref{*}, 
we conclude that 
\begin{align*}
\sum_{i=1}^r
\Bigl(\|\Phi_i\|^2+&
n_iH_i^2\bigl(1-\frac{n_i}{n-2}\bigr)\Bigr)
+
\frac{2(n-3)}{n-2}
\sum_{1\leq i<j\leq r}
\sum_{k=1}^{n_i}\sum_{\ell=1}^{n_j}
\Bigl\|
\alpha_f(\sigma_{i*}e_k^i,
\sigma_{j*}e_\ell^j)
\Bigr\|^2
\\
&\leq
4c-\frac{2c}{n-2}
\sum_{1\leq i<j\leq r}n_in_j
=
\frac{c}{n-2}
\Bigl(\sum_{i=1}^r n_i^2-n^2+4n-8\Bigr),
\end{align*}
which proves \eqref{pro}. 
\qed

\begin{lemma}\label{NR}
Let $f\colon M^n\to\bar {M}^{n+m}$, 
$n\geq4$, be an isometric immersion into a 
Riemannian manifold whose curvature 
operator is bounded below by a constant 
$c$. Then the following hold: 
\begin{enumerate}[
label=(\roman*),
topsep=1pt,
itemsep=0.2ex,
leftmargin=*,
partopsep=0.2ex,
parsep=0.2ex,
]

\item Suppose that $\delta$ is a Dupin 
principal normal vector field of $f$ and that 
$TM$ admits an orthogonal splitting 
$TM=D_1\oplus D_2$, where $D_1$ is 
a rank-two subbundle satisfying 
$D_1\subset\mathcal D_f$ 
and $D_2\subset E_\delta$. Then the 
curvature operator of $M^n$ is also 
bounded below by $c$. 
\item Suppose that $n=4$, and at 
each point $x\in M^4$ where $f$ is 
not totally umbilical, there exist two 
distinct Dupin principal normals 
$\delta_1$ and $\delta_2$, 
both of multiplicity $2$, such that 
$$
\<\delta_1,\delta_2\>(x)=-c
\quad \text{ and }\quad 
T_xM=E_{\delta_1}(x)\oplus E_{\delta_2}(x).
$$
If $c\ge0$, then $M^4$ has nonnegative 
curvature operator. 
\end{enumerate}
\end{lemma}
\proof
\noindent $(i)$ Fix $x\in M^n$ and choose an 
orthonormal basis $\{e_i\}_{i=1}^n$ of $T_xM$ 
such that $D_1(x)=\spa\{e_1,e_2\}$ 
and $D_2(x)=\spa\{e_3,\dots,e_n\}$. 
By our assumptions, 
\[
\alpha_f(e_i,e_j)=
\begin{cases}
\delta(x),& i=j\geq3,\\
0,&\text{otherwise}.
\end{cases}
\]

Let 
\[
V=\sum_{1\leq i<j\leq n}
v_{ij}e_i\wedge e_j
\]
be an arbitrary element of $\Lambda^2T_xM$. 
Using \eqref{RR}, we obtain 
\begin{align*}
\langle\!\langle\mathcal R_MV,V\rangle\!\rangle
=
\langle\!\langle
{\mathcal R}_fV,V
\rangle\!\rangle
+\|\delta(x)\|^2
\sum_{3\leq i<j\leq n}v_{ij}^2\geq c\|V\|^2,
\end{align*}
since the curvature operator of the ambient 
manifold is bounded below by $c$. 
Thus, the curvature operator of $M^n$ 
is also bounded below by $c$. 
\vspace{1ex}

\noindent $(ii)$ Fix $x\in M^4$. 
Suppose first that $f$ is not totally 
umbilical at $x$. Choose an orthonormal 
basis $\{e_i\}_{i=1}^4$ of $T_xM$ such that 
$E_{\delta_1}(x)=\spa\{e_1,e_2\}$ and 
$E_{\delta_2}(x)=\spa\{e_3,e_4\}$. 
Let 
\[
V=\sum_{1\leq i<j\leq4}
v_{ij}e_i\wedge e_j
\]
be an arbitrary element of $\Lambda^2T_xM$. 
Using \eqref{RR} and our assumptions, 
we obtain 
\begin{align*}
\langle\!\langle\mathcal R_MV,V\rangle\!\rangle
&=
\langle\!\langle
{\mathcal R}_fV,V
\rangle\!\rangle
+v_{12}^2\|\delta_1\|^2+v_{34}^2\|\delta_2\|^2+
\langle\delta_1,\delta_2\rangle\sum_{i=1}^2\sum_{j=3}^4
v_{ij}^2\\
&\geq
c\|V\|^2
+v_{12}^2\|\delta_1\|^2+v_{34}^2\|\delta_2\|^2
-c\sum_{i=1}^2\sum_{j=3}^4v_{ij}^2
\\
&=
v_{12}^2(c+\|\delta_1\|^2)
+v_{34}^2(c+\|\delta_2\|^2)\ge0,
\end{align*}
where the last inequality follows from 
$c\geq0$. 

Suppose now that $f$ is totally umbilical 
at $x$. Then 
\[
\alpha_f(X,Y)
=\langle X,Y\rangle\mathcal H_f(x)
\]
for all $X,Y\in T_xM$. Hence, by 
the Gauss equation, 
\[
\langle\!\langle\mathcal R_MV,V\rangle\!\rangle
=
\langle\!\langle
{\mathcal R}_fV,V
\rangle\!\rangle
+H^2(x)\|V\|^2
\geq c\|V\|^2\geq0.
\]
Therefore, $M^4$ has nonnegative 
curvature operator. 
\qed
\vspace{1ex}

The next result describes the extrinsic geometry 
of isometric immersions of Riemannian products 
satisfying the pinching condition \eqref{1}, 
under suitable restrictions on the dimensions 
of the factors. 

\begin{proposition}\label{propprod}
Let $f\colon M^n\to\bar {M}^{n+m}$, 
$n\geq4$, be an isometric immersion of a 
Riemannian manifold $M^n$ that is isometric 
to a Riemannian product 
\[
(N_1^{n_1},g_1)\times\cdots\times
(N_r^{n_r},g_r),
\qquad r\geq2,
\]
where $1\leq n_i\leq n-2$ for all 
$1\leq i\leq r$, and at least one factor 
has dimension at least two. 
Suppose that the curvature operator of 
$\bar {M}^{n+m}$ is bounded below 
by a constant $c$. If inequality \eqref{1} 
holds at every point, then the following 
assertions hold: 
\begin{enumerate}[
label=(\roman*),
topsep=1pt,
itemsep=0.2ex,
leftmargin=*,
partopsep=0.2ex,
parsep=0.2ex,
]

\item If $n_i<n-2$ for all 
$1\leq i\leq r$, then $c\leq0$. 
Moreover, if $c\geq0$, then necessarily 
$c=0$, the immersion $f$ is totally 
geodesic, $M^n$ has nonnegative curvature 
operator, and equality holds in \eqref{1} 
at every point. 

\item If $n\geq5$, $r=2$, $n_1=2$, 
and $n_2=n-2$, then $c\leq0$, 
$TN_1\subset\mathcal D_f$, and 
$\delta=n\mathcal H_f/(n-2)$ is a 
Dupin principal normal of $f$ 
satisfying $TN_2\subset E_\delta$. 
Moreover, if $c\geq0$, then necessarily 
$c=0$, $M^n$ has nonnegative curvature 
operator, and equality holds in \eqref{1} 
at every point. 

\item If $n=4$, $r=2$, and 
$n_1=n_2=2$, then there exist smooth 
Dupin principal normal vector fields 
$\delta_1$ and $\delta_2$ of $f$ such that 
\[
TN_i\subset E_{\delta_i},\quad i=1,2, \qquad
\langle\delta_1,\delta_2\rangle=-c.
\]
Furthermore, equality holds in \eqref{1} 
at every point. If $c\geq0$, then 
$M^4$ has nonnegative curvature operator. 
\end{enumerate}
\end{proposition}
\proof
Suppose that $f$ satisfies inequality \eqref{1} 
at every point. By our dimension assumptions, 
the left-hand side of \eqref{pro} is 
nonnegative, and consequently, 
\be\label{min}
c\Bigl(\sum_{i=1}^r n_i^2-n^2+4n-8\Bigr)\geq0.
\ee

\noindent $(i)$ Suppose that $n_i<n-2$ for all 
$1\leq i\leq r$. Choose a factor of dimension 
at least two, say $n_1\geq2$, 
and write $n_1=n-k$, where $3\leq k\leq n-2$. 
Then 
\[
\sum_{i=1}^r n_i^2
=(n-k)^2+\sum_{i=2}^r n_i^2
\leq(n-k)^2+k^2
\leq(n-2)^2+4.
\]
The last inequality follows from the convexity 
of the function $t\mapsto(n-t)^2+t^2$ 
on the interval $[3,n-2]$. Moreover, equality 
cannot hold throughout, since equality in the 
first inequality requires $r=2$, while equality 
in the second requires $k=n-2$, which would 
imply $n_2=n-2$, contrary to our assumption. 
Thus, 
\[
\sum_{i=1}^r n_i^2<n^2-4n+8,
\]
and \eqref{min} yields $c\leq0$. 

Suppose now that $c\geq0$. Then 
necessarily $c=0$. Since $n_i<n-2$ 
for all $i$, it follows from \eqref{pro} that 
$\Phi_i=H_i=0$ for every $1\leq i\leq r$, 
and
\[
\alpha_f(TN_i,TN_j)=0
\quad\text{for all }i\neq j.
\]
Consequently, each immersion $f_i$ is totally 
geodesic, and the second fundamental form 
of $f$ is adapted to the product net 
$\{TN_i\}_{i=1}^r$. By \eqref{ai}, the 
immersion $f$ is therefore totally geodesic. 
The Gauss equation then implies that $M^n$ 
has nonnegative curvature operator. Moreover, 
since $S=H=c=0$, equality holds in \eqref{1} 
everywhere. 
\vspace{1ex}

\noindent $(ii)$ Suppose that $n\geq5$, 
$n_1=2$, and $n_2=n-2$. Since the 
right-hand side of \eqref{pro} 
vanishes and all terms on its left-hand side
are nonnegative, each of these terms vanishes. 
Hence $\Phi_1=\Phi_2=0$, $H_1=0$, 
and $\alpha_f(TN_1,TN_2)=0$. 
Consequently, $f_1$ is totally geodesic, 
$f_2$ is totally umbilical, and the second 
fundamental form of $f$ is adapted to the 
product net $\{TN_1,TN_2\}$. 
It then follows from \eqref{ai} that 
$TN_1\subset\mathcal D_f$, while 
$\delta=n\mathcal H_f/(n-2)$ is a 
Dupin principal normal of $f$ satisfying 
$TN_2\subset E_\delta$. 

Now let $x\in M^n$, and choose unit 
vectors $v_i\in T_xN_i$, $i=1,2$. 
Since $M^n$ is a Riemannian product, 
the mixed sectional curvature vanishes. 
Moreover, $\a_f(v_1,v_1)=0$ and 
$\a_f(v_1,v_2)=0$. Thus, the Gauss 
equation gives 
\[
0=K_M(v_1\wedge v_2)
=K_{\bar M}
(f_*v_1\wedge f_*v_2)\geq c.
\]
Therefore, $c\leq0$. 

Suppose now that $c\geq0$. Then $c=0$. 
Furthermore, since $S=(n-2)\|\delta\|^2$, 
and $n\mathcal H_f=(n-2)\delta$, 
we obtain 
\[
S=\frac{n^2H^2}{n-2}.
\]
Hence, equality holds in \eqref{1} 
at every point. Finally, part $(i)$ of 
Lemma \ref{NR} implies that $M^n$ 
has nonnegative curvature operator. 
\vspace{1ex}

\noindent $(iii)$ Suppose that $n=4$ and 
$n_1=n_2=2$. Then the right-hand side 
of \eqref{pro} vanishes, yielding 
$\Phi_1=\Phi_2=0$ and $\a_f(TN_1,TN_2)=0$.
Thus, both immersions $f_1$ and $f_2$ 
are totally umbilical, and the second 
fundamental form of $f$ is adapted to 
the product net $\{TN_1,TN_2\}$. 
Consequently, by \eqref{ai}, the smooth 
normal vector fields $\delta_i=\mathcal H_{f_i}$, 
$i=1,2$, are Dupin principal normals of $f$ 
satisfying $TN_i\subset E_{\delta_i}$. 

Let $x\in M^4$, and choose unit vectors 
$v_i\in T_xN_i$, $i=1,2$. Since the mixed 
sectional curvature of a Riemannian 
product vanishes, the Gauss equation yields 
\[
0=K_M(v_1\wedge v_2)
=K_{\bar M}
(f_*v_1\wedge f_*v_2)
+\langle\delta_1,\delta_2\rangle
\geq c+\langle\delta_1,\delta_2\rangle.
\]
Therefore, 
$\langle\delta_1,\delta_2\rangle\leq-c$.
On the other hand, since 
$$
S=2\left(\|\delta_1\|^2+\|\delta_2\|^2\right)
\quad\text{ and }\quad
4H^2=\|\delta_1\|^2+\|\delta_2\|^2
+2\langle\delta_1,\delta_2\rangle,
$$
inequality \eqref{1}, which in this 
case reads $S\leq4c+8H^2$, implies 
$\langle\delta_1,\delta_2\rangle\geq-c$. 
Consequently, 
$\langle\delta_1,\delta_2\rangle=-c$, 
and equality holds in \eqref{1} 
at every point. 

If $c\geq0$, then part $(ii)$ of Lemma 
\ref{NR} implies that $M^4$ has 
nonnegative curvature operator. 
This completes the proof.
\qed

%\vspace{-5pt}
\subsection{The Bochner--Weitzenb\"ock operator}

Let $(M^n,\<\cdot,\cdot\>)$ be a 
Riemannian manifold of dimension $n\geq4$, 
with Levi-Civit\'a connection $\n$ and 
curvature tensor $R$. The Ricci tensor is 
defined by 
\[
\Ric(X,Y)=
\sum_{i=1}^n\<R(X,E_i)E_i,Y\>,
\qquad X,Y\in\mathcal X(M),
\]
where $\{E_i\}_{i=1}^n$ is a local 
orthonormal frame. We also denote by 
$\Ric$ the corresponding Ricci endomorphism. 

At each point $x\in M^n$, the 
{\emph{Bochner--Weitzenb\"ock operator}} 
$\B^{[2]}$ is the self-adjoint 
endomorphism of $\Lambda^2T_xM$ defined by 
the following formula (see \cite{sea}): 
\begin{align}\label{WB}
\<\<\B^{[2]}(v_1\wedge v_2), w_1\wedge w_2\>\>
=&\;\Ric(v_1,w_1)\<v_2,w_2\>+\Ric(v_2,w_2)\<v_1,w_1\>\nonumber\\
&-\Ric(v_1,w_2)\<v_2,w_1\>-\Ric(v_2,w_1)\<v_1,w_2\>\nonumber\\
&-2\<R(v_1,v_2)w_2,w_1\>,
\end{align}
and 
extended linearly to all of $\Lambda^2T_xM$. 
Here $\<\<\cdot,\cdot\>\>$ denotes the 
inner product on $\Lambda^2T_xM$ induced 
by the Riemannian metric, given by 
\[
\<\<v_1\wedge v_2,w_1\wedge w_2\>\>
=\det(\<v_i,w_j\>).
\]

Each $Z\in\Lambda^2T_xM$ determines 
a 2-form $\omega$ by 
\[
\omega(v,w)=\<\<Z,v\wedge w\>\>.
\]
Equivalently, $Z$ can be identified 
with the skew-symmetric endomorphism 
of $T_xM$ defined by 
\[
\<Z(v),w\>=\<\<Z,v\wedge w\>\>.
\]

Under the natural identification of 
2-vectors with 2-forms, $\B^{[2]}$ 
also acts on the space $\Omega^2(M)$ 
of smooth 2-forms. More precisely, 
\begin{align*}
\B^{[2]}\omega(X_1,X_2)
=\omega(\Ric(X_1),X_2)
+\omega(X_1,\Ric(X_2))-\sum_{i=1}^n
\omega(R(X_1,X_2)E_i,E_i).
\end{align*}
These two descriptions are related by 
\[
\omega(\B^{[2]}(X_1\wedge X_2))
=(\B^{[2]}\omega)(X_1,X_2).
\]
In particular, taking $\omega$ to be 
the 2-form dual to $w_1\wedge w_2$ 
recovers \eqref{WB}. Throughout the paper, 
we will generally identify 2-forms 
with 2-vectors. 

Recall that the Bochner--Weitzenb\"ock 
formula for 2-forms is 
\[
\<\Delta\omega,\omega\>
=\frac12\Delta\|\omega\|^2
+\|\n\omega\|^2
+\<\B^{[2]}\omega,\omega\>,
\]
where $\Delta=d\delta+\delta d$ 
denotes the Hodge Laplacian, with 
the corresponding sign convention 
for the scalar Laplacian. 
Consequently, if $M$ is compact 
and $\B^{[2]}$ is nonnegative, 
then every harmonic 2-form is parallel. 

In dimension four, nonnegative isotropic 
curvature at a point is equivalent to the 
positive semidefiniteness of the 
Bochner--Weitzenb\"ock 
operator $\B^{[2]}$ at that point 
(see, for instance, \cite{MWolf}).

\subsection{Even-dimensional submanifolds} 
It is known that even-dimensional Riemannian 
manifolds with nonnegative isotropic curvature 
have positive semidefinite Bochner--Weitzenb\"ock 
operator $\B^{[2]}$ (see \cite{sea}). 
In this section, we strengthen this 
result in a form adapted to our purposes. 

Let $(M^n,g)$, $n=2k$, be an even-dimensional 
Riemannian manifold, and let 
$\omega\in\Omega^2(M^n)$ be a $2$-form. 
At each point $x\in M^n$, there exists 
an orthonormal basis 
$\{v_1,w_1,\ldots,v_k,w_k\}$ of $T_xM$
such that 
\be\label{Lambda}
\omega(v_a,w_b)=\lambda_a\delta_{ab},
\qquad
\omega(v_a,v_b)=\omega(w_a,w_b)=0,
\quad 1\leq a,b\leq k,
\ee
where $\lambda_1,\ldots,\lambda_k$ are 
nonnegative real numbers. 

We recall the following lemma from 
\cite{v5}.

\begin{lemma}\label{lemmaeven2} 
Let $M^n$, with $n=2k\geq4$, be an 
even-dimensional Riemannian manifold. 
Suppose that $M^n$ has nonnegative 
isotropic curvature at a point $x\in M^n$. 
Then the following assertions hold: 
\begin{enumerate}[
label=(\roman*),
topsep=1pt,
itemsep=0.2ex,
leftmargin=*,
partopsep=0.2ex,
parsep=0.2ex,
]

\item The Bochner--Weitzenb\"ock 
operator $\B^{[2]}$ is positive 
semidefinite at $x$. 

\item Let $\omega\in\Omega^2(M^n)$ 
be a 2-form satisfying 
$\<\B^{[2]}\omega,\omega\>(x)=0$. 
Choose an orthonormal basis
$\{v_1,w_1,\dots,v_k,w_k\}$ 
of $T_xM$ such that \eqref{Lambda} holds. 
Then, for every $1\leq a,b\leq k$ 
with $a\neq b$, we have 
\begin{align}\label{even=}
\lambda_a^2\Big(\left<R(v_b,v_a)v_a,v_b\right>&
+\left<R(v_b,w_a)w_a,v_b\right> +\left<R(w_b,v_a)v_a,w_b\right> \nonumber\\
&+\left<R(w_b,w_a)w_a,w_b\right> -2\left|\<R(v_b,w_b)v_a,w_a\>\right|\Big)=0.
\end{align} 
\end{enumerate}
\end{lemma}

\begin{proposition}\label{propeven}
Let $M^n$, with $n=2k\geq4$, be an even-dimensional 
Riemannian manifold. Suppose that an isometric immersion 
$f\colon M^n \to\bar {M}^{n+m}$ 
satisfies inequality \eqref{1} at a point $x\in M^n$, where 
$c$ is a lower bound on the curvature operator of 
$\bar {M}^{n+m}$. 
Then the following assertions hold: 
\begin{enumerate}[
label=(\roman*),
topsep=1pt,
itemsep=0.2ex,
leftmargin=*,
partopsep=0.2ex,
parsep=0.2ex,
]

\item The Bochner--Weitzenb\"ock operator 
$\B^{[2]}$ is positive semidefinite at $x$. 
\item Let $\omega\in\Omega^2(M^n)$ 
be a 2-form such that 
$\<\B^{[2]}\omega,\omega\>(x)=0$ and $\omega_x\neq0$. 
If $n\geq6$, then the following conditions are 
satisfied at $x$: 
\begin{enumerate}[
label=(\alph*),
topsep=1pt,
itemsep=0.2ex,
leftmargin=*,
partopsep=0.2ex,
parsep=0.2ex,
]

\item The inequality \eqref{1} holds as an 
equality at $x$, and $c=0$. 
\item If $f$ is not totally geodesic at $x$, 
then exactly one of the numbers $\lambda_a$, 
$1\leq a\leq k$, in \eqref{Lambda} is nonzero. 
Moreover, the normal vector
$\delta=n\mathcal H_f/(n-2)$ is a Dupin 
principal normal of $f$ at $x$ of multiplicity 
$n-2$, and the relative nullity space is given 
by $\mathcal D_f(x)=E^\perp_\delta(x)$. 
\end{enumerate}
\end{enumerate}
\end{proposition}
\proof
\noindent $(i)$ This follows directly from part 
$(i)$ of Proposition \ref{propu} together with 
part $(i)$ of Lemma \ref{lemmaeven2}. 
\vspace{1ex}

\noindent $(ii)$ Choose an orthonormal basis 
$\{e_i\}_{i=1}^n$ of 
$T_xM$ such that $e_{2a-1}=v_a$ and 
$e_{2a}=w_a$ for $1\leq a\leq k$, where 
$\{v_1,w_1,\dots,v_k,w_k\}$ is chosen so that 
\eqref{Lambda} holds. Since $\omega_x\neq0$, 
we may assume, without loss of generality, 
that $\lambda_1\neq0$. 
It follows from \eqref{even=} that, for every
$a\neq1$, there exists $\varepsilon_a=\pm1$
such that the orthonormal four-frame
$\{v_1,w_1,v_a,w_a\}$ satisfies
\eqref{eq1} with $\varepsilon=\varepsilon_a$.

%It follows from 
%\eqref{even=} that the isotropic curvature 
%vanishes for the four-frame 
%$\{v_1,w_1,v_a,w_a\}$ for every $a\neq1$. 

By part $(ii)$ of Proposition \ref{propu}, inequality 
\eqref{1} holds as an equality at $x$, and the normal vector 
$\delta=n\mathcal H_f/(n-2)$ is a Dupin principal 
normal of $f$ such that $E_\delta(x)$ contains 
the orthogonal complement of 
${\operatorname{span}}\{v_1,w_1,v_a,w_a\}$ 
in $T_xM$ for every $a\neq1$. Moreover, conditions 
\eqref{c1}--\eqref{c2} are satisfied for all such four-frames. 
Since $k\geq3$, every vector $e_i$, 
$3\leq i\leq n$, belongs to one of these 
orthogonal complements. Consequently, 
$
{\operatorname{span}}\left\{e_i\right\}_{i=3}^n\subset E_\delta(x). 
$
In particular, 
$\a_{ii}=\delta$ for all $i\geq3$. 
It follows from the first equality in \eqref{c1} and 
\eqref{c2} that 
$
\a_{11}=\a_{22}=0.
$
Furthermore, the inclusion above, 
together with the last three equalities 
in \eqref{c1}, yields 
$$
\a_{ij}=0\quad\text{ for all}\quad1\leq i\neq j\leq n.
$$
Hence, ${\operatorname{span}}\{e_1,e_2\}\subset\mathcal D_f(x)$, 
and therefore 
$$
S(x)=(n-2)\|\delta\|^2=\frac{n^2H^2(x)}{n-2}.
$$ 
Since 
equality holds in \eqref{1}, we conclude that 
$c=0$. This proves $(a)$. 

Suppose now that $f$ is not totally geodesic at $x$. 
Then $\delta\neq0$ is a Dupin 
principal normal of multiplicity $n-2$ with 
\be\label{Ed}
E_\delta(x)={\operatorname{span}}\{e_3,\dots,e_n\},
\ee
and the relative nullity subspace of $f$ at $x$ is 
$\mathcal D_f(x)={\operatorname{span}}\{e_1,e_2\}$. 

It remains to prove that exactly one 
of the numbers $\lambda_a$, 
$1\leq a\leq k$, is nonzero. Suppose, to the 
contrary, that $\lambda_2\neq0$ as well. Then 
\eqref{even=} implies that the isotropic 
curvature vanishes for the four-frame 
$\{v_2,w_2,v_a,w_a\}$ for every $a\neq2$. 
Applying again part $(ii)$ of Proposition \ref{propu}, 
we conclude that $E_\delta(x)$ contains the 
orthogonal complement of 
${\operatorname{span}}\{v_2,w_2,v_a,w_a\}$ 
for every $a\neq2$. Since $k\geq3$, 
we may choose $a\geq3$, obtaining 
${\operatorname{span}}
\left\{e_1,e_2\right\}\subset E_\delta(x)$. 
This contradicts \eqref{Ed}. Therefore, $\lambda_a=0$ 
for all $2\leq a\leq k$, 
which completes the proof. 
\qed

\section{Proof of the main result for the case $n=4$}

\subsection{Geometry of four-dimensional manifolds}
In this subsection, we collect some basic facts about 
four-dimensional geometry. For a detailed exposition 
of the subject, we refer the reader to \cite{LeB}. 

The space of smooth 2-forms on an oriented 
four-dimensional Riemannian manifold 
$(M,\<\cdot,\cdot\>)$ decomposes as a direct sum 
\[
\Omega^2(M)=\Omega_+^2(M)\oplus\Omega_-^2(M)
\]
of the eigenspaces of the Hodge star operator 
$\ast\colon\Omega^2(M)\to\Omega^2(M)$. 
The elements of $\Omega_+^2(M)$ are called 
\emph{self-dual 2-forms}, whereas those of 
$\Omega_-^2(M)$ are called 
\emph{anti-self-dual 2-forms}. 
Accordingly, at every point $x\in M$, we have 
the splitting 
$$
\Lambda^2T_xM=\Lambda_+^2T_xM\oplus\Lambda_-^2T_xM, 
$$
where $\Lambda_\pm^2T_xM$ 
are the eigenspaces of the Hodge star operator 
$\ast\colon\Lambda^2T_xM\to\Lambda^2T_xM$. 
Both spaces $\Lambda_+^2T_xM$ and $\Lambda_-^2T_xM$ 
are $\B^{[2]}$-invariant (see Proposition 1 in \cite{smichig}). 
Hence, we obtain the corresponding decomposition 
$\B^{[2]}=\B_+^{[2]}\oplus\B_-^{[2]}$. 

With respect to the splitting 
$\Lambda^2TM=\Lambda^2_+TM\oplus\Lambda^2_-TM$, 
the curvature operator of $M$ is 
$$
\mathcal R=
\begin{pmatrix}
W^++\dfrac{\tau}{12}I&B\\[2mm]
B^*&W^-+\dfrac{\tau}{12}I
\end{pmatrix},
$$
where $\tau$ is the scalar curvature, 
$$
W^\pm+\dfrac{\tau}{12}I\colon\Lambda^2_\pm TM\to\Lambda^2_\pm TM,
$$
$W^+$ (respectively, $W^-$) is the trace-free block called 
\emph{self-dual Weyl operator} (respectively, anti-self-dual 
Weyl operator), and $B$ is the trace-free Ricci tensor. 
Moreover, 
$$
\B_\pm^{[2]}=\dfrac{\tau}{3}I-2W^\pm. 
$$

Suppose now that $M$ is a compact, 
oriented Riemannian four-manifold. The Hodge 
theorem guarantees that every de Rham 
cohomology class on $M$ has a unique 
harmonic representative. In particular, the space 
$\mathscr H^2(M)$ of harmonic 2-forms 
decomposes as 
$$
\mathscr H^2(M)=\mathscr H^2_+(M)
\oplus\mathscr H^2_-(M),
$$ 
where $\mathscr H^2_+(M)$ and 
$\mathscr H^2_-(M)$ denote the 
spaces of self-dual and anti-self-dual 
harmonic 2-forms, respectively. The 
dimensions 
$b_\pm(M)=\dim\mathscr H^2_\pm(M)$ 
are oriented homotopy invariants of $M$. 
Their difference 
$
\sigma(M)=b_+(M)-b_-(M)
$ 
is the signature of $M$, 
while their sum equals the second Betti number 
$b_2(M)$ of the manifold. 

\subsection{Auxiliary results for four-dimensional submanifolds}

We now state the following auxiliary results. 

\begin{lemma}\label{proalg}
Let $(V^4,\langle\cdot,\cdot\rangle)$ be an oriented 
four-dimensional Euclidean vector space, and let 
$\omega\in\Lambda^2_+V^*$ be a self-dual 
two-form of length $\|\omega\|=\sqrt{2}$. 
Then the following assertions hold: 
\begin{enumerate}[
label=(\roman*),
topsep=1pt,
itemsep=0.2ex,
leftmargin=*,
partopsep=0.2ex,
parsep=0.2ex,
]

\item There 
exists a positively oriented orthonormal basis 
$\{e_i\}_{i=1}^4$ of $V$ such that 
$$
\omega=\omega_1\wedge \omega_2
+\omega_3\wedge \omega_4,
$$
where $\{\omega_i\}_{i=1}^4$ is the dual basis. 
\item The skew-symmetric endomorphism 
associated with $\omega$ is an orthogonal complex 
structure $J$ on $V$, compatible with its orientation, and
$$
\langle Jv,w\rangle=\omega(v,w),\qquad v,w\in V.
$$
\end{enumerate}
\end{lemma}
\proof
\noindent $(i)$ Associate to $\omega$ the 
skew-symmetric endomorphism $A\colon V\to V$ 
defined by 
\[
\langle Av,w\rangle=\omega(v,w),
\qquad v,w\in V.
\]

Choose any unit vector $e_1\in V$ and extend 
it to a positively oriented orthonormal basis 
$\{e_1,v_2,v_3,v_4\}$ of $V$, with dual basis 
$\{\phi_1,\phi_2,\phi_3,\phi_4\}$. 
Since $\omega$ is self-dual, we may write 
$$
\omega=
x_1(\phi_1\wedge\phi_2+
\phi_3\wedge\phi_4)
+x_2(\phi_1\wedge\phi_3-
\phi_2\wedge\phi_4)
+x_3(\phi_1\wedge\phi_4+
\phi_2\wedge\phi_3).
$$
The three self-dual two-forms appearing 
above are mutually orthogonal and each 
has squared norm $2$. Since 
$\|\omega\|^2=2$, it follows that 
$x_1^2+x_2^2+x_3^2=1$. 

By the definition of $A$, 
\[
Ae_1=x_1v_2+x_2v_3+x_3v_4.
\]
Consequently, $\|Ae_1\|=1$. 
Setting $e_2=Ae_1$, we obtain 
$\langle e_1,e_2\rangle=0$, since $A$ 
is skew-symmetric. Thus, $\{e_1,e_2\}$ 
is an orthonormal 
set. Moreover, 
$\omega(e_1,e_2)
=\langle Ae_1,e_2\rangle=1$. 

Complete $\{e_1,e_2\}$ to a positively 
oriented orthonormal basis 
$\{e_i\}_{i=1}^4$ of $V$, and let 
$\{\omega_i\}_{i=1}^4$ 
be its dual basis. 
By self-duality, we may write 
$$
\omega=
y_1(\omega_1\wedge\omega_2+
\omega_3\wedge\omega_4)
+y_2(\omega_1\wedge\omega_3-
\omega_2\wedge\omega_4)
+y_3(\omega_1\wedge\omega_4+
\omega_2\wedge\omega_3).
$$
Now $y_1=\omega(e_1,e_2)=1$, and the
normalization $\|\omega\|^2=2$ gives
$y_1^2+y_2^2+y_3^2=1$. 
Consequently, $y_2=y_3=0$, and therefore
\[
\omega=\omega_1\wedge\omega_2+
\omega_3\wedge\omega_4,
\]
which proves $(i)$. 
\vspace{1ex}

\noindent $(ii)$ By part $(i)$ and the 
definition of $A$, the positively oriented 
orthonormal basis constructed above 
satisfies 
$$
Ae_1=e_2,\quad Ae_2=-e_1,\quad 
Ae_3=e_4,\quad Ae_4=-e_3.
$$
Hence, $J=A$ is an orthogonal 
complex structure on $V$ satisfying 
\[
\langle Jv,w\rangle=\omega(v,w),
\qquad v,w\in V.
\]
Moreover, the ordered basis 
$\{e_1,Je_1,e_3,Je_3\}
=\{e_1,e_2,e_3,e_4\}$ 
is positively oriented. Therefore, 
$J$ induces the given orientation 
on $V$, completing the proof. 
\qed

\paragraph{Remark (anti-self-dual case).}
If $\omega\in\Lambda^2_-V^*$ satisfies 
$\|\omega\|=\sqrt{2}$, the same argument, 
with the appropriate signs changed, 
gives a positively oriented orthonormal 
basis in which 
\[
\omega=\omega_1\wedge\omega_2
-\omega_3\wedge\omega_4.
\]
The associated skew-symmetric endomorphism 
$A$ is again an orthogonal complex structure, 
but it induces the orientation opposite 
to the given one.

\begin{proposition}\label{beta}
Let $(M^4,g)$ be a compact, oriented 
Riemannian four-manifold with nonnegative isotropic 
curvature. If $b_+(M)>0$, (respectively, $b_-(M)>0$), 
then at every $x\in M^4$ there exists a positively 
oriented orthonormal basis 
$\{e_i\}_{i=1}^4$ of $T_xM$ such that 
\begin{equation}\label{eq:zero-ic}
R_{1313}+R_{1414}+R_{2323}+R_{2424}
-2\varepsilon R_{1234}=0,
\end{equation}
where $\varepsilon=+1$ (respectively, 
$\varepsilon=-1$). 
\end{proposition}
\proof
We first consider the case $b_+(M)>0$; 
the case $b_-(M)>0$ is analogous. 
By Hodge theory, there exists a nonzero 
harmonic self-dual $2$-form $\omega$, that is, 
$\Delta\omega=0$ and $*\omega=\omega$. 
The Bochner--Weitzenb\"ock formula 
for self-dual $2$-forms gives 
$$
\Delta\omega=\nabla^*\nabla\omega+
\left(\frac{\tau}{3}I-2W^+\right)\omega,
$$
where $\tau$ is the scalar curvature 
and $W^+$ is the self-dual Weyl 
curvature operator. 

Since $M^4$ has nonnegative isotropic 
curvature, the symmetric operators 
\[
\frac{\tau}{3}I-2W^\pm\colon
\Lambda^2_\pm T_x^*M
\longrightarrow\Lambda^2_\pm T_x^*M
\]
are positive semidefinite at every point $x\in M^4$. 

Integrating the Bochner--Weitzenb\"ock formula 
yields 
$$
0=\int_M\left(\|\nabla\omega\|^2+
\left\langle\left(\frac{\tau}{3}I-2W^+\right)\omega,
\omega\right\rangle\right)\,dV_g.
$$
Consequently, 
\begin{equation}\label{eq:parallel-kernel}
\nabla\omega=0,
\qquad
\left(\frac{\tau}{3}I-2W^+\right)\omega=0.
\end{equation}
Since $\omega$ is nonzero and parallel, 
its norm is a positive constant. After 
rescaling $\omega$, we may assume that 
$\|\omega\|=\sqrt{2}$ everywhere. 
Fix $x\in M^4$. By Lemma \ref{proalg}, 
there exists a positively oriented 
orthonormal basis 
$\{e_i\}_{i=1}^4$ of $T_xM$ 
such that 
\[
\omega_x=\omega_1\wedge\omega_2
+\omega_3\wedge\omega_4,
\]
where $\{\omega_i\}_{i=1}^4$ 
is the dual basis. 

Identifying $2$-forms with $2$-vectors 
via the metric, the bivector 
$\xi=e_1\wedge e_2+e_3\wedge e_4$
satisfies, by \eqref{eq:parallel-kernel}, 
\[
\left(\frac{\tau}{3}I-2W^+\right)\xi=0.
\]
The four-dimensional curvature identity 
then gives 
\begin{align*}
&R_{1313}+R_{1414}+R_{2323}+R_{2424}
-2R_{1234}=
\frac12\left\langle\!\left\langle
\left(\frac{\tau}{3}I-2W^+\right)\xi,
\xi
\right\rangle\!\right\rangle=0.
\end{align*}
Thus, \eqref{eq:zero-ic} holds at every 
point $x\in M^4$ with $\varepsilon=+1$. 

If $b_-(M)>0$, the same argument 
shows that there exists a nonzero 
parallel anti-self-dual $2$-form 
$\omega$ satisfying 
\[
\left(\frac{\tau}{3}I-2W^-\right)\omega=0.
\]
After normalization, at every point 
$x\in M^4$ we may choose a positively 
oriented orthonormal basis such that 
\[
\omega_x=\omega_1\wedge\omega_2
-\omega_3\wedge\omega_4.
\]
Setting $\xi=e_1\wedge e_2-e_3\wedge e_4$, 
we obtain 
\[
\left(\frac{\tau}{3}I-2W^-\right)\xi=0.
\]
Consequently, 
\begin{align*}
&R_{1313}+R_{1414}+R_{2323}+R_{2424}
+2R_{1234}=
\frac12\left\langle\!\left\langle
\left(\frac{\tau}{3}I-2W^-\right)\xi,
\xi
\right\rangle\!\right\rangle=0.
\end{align*}
This completes the proof. 
\qed
\vspace{1ex}

The following result will be used in the sequel.

\begin{theorem}\label{n=4beta}
Let $f\colon M^4\to\bar {M}^{4+m}$ 
be an isometric immersion of a compact, 
oriented Riemannian four-manifold 
satisfying inequality \eqref{1} at every point, 
where $c\geq0$ is a lower bound on the 
curvature operator of $\bar {M}^{4+m}$. 
If $b_2(M)>0$, then equality holds 
in \eqref{1} at every point, and $M^4$ 
admits a K\"ahler structure with respect 
to which $f$ is $(2,0)$-geodesic. 

Moreover, if the scalar curvature of $M^4$ 
vanishes identically, then $c=0$, the 
immersion $f$ is totally geodesic, and 
$M^4$ is flat, with $2\leq b_1(M)\leq4$. 
In particular, $b_1(M)=4$ if and only 
if $M^4$ is isometric to a flat torus. 
\end{theorem}
\begin{proof}
By Propositions \ref{propu} and \ref{propeven}, 
the manifold $M^4$ has nonnegative isotropic 
curvature, and the Bochner--Weitzenb\"ock operator 
$\B^{[2]}$ is positive semidefinite. Since $b_2(M)>0$, 
either $b_+(M)>0$ or $b_-(M)>0$. 
Moreover, Proposition \ref{beta} guarantees 
the existence of an orthonormal four-frame 
with vanishing isotropic curvature at every 
point of $M^4$. Hence, part $(ii)$ of 
Proposition \ref{propu} implies that equality 
holds identically in \eqref{1}. 

We first consider the case $b_+(M)>0$. 
By Hodge theory, there exists a nonzero 
harmonic self-dual $2$-form $\omega$. 
The Bochner--Weitzenb\"ock formula gives 
\[
0=\int_M\left(
\|\nabla\omega\|^2+
\langle\B^{[2]}\omega,\omega\rangle
\right)\,dV_g.
\]
Since $\B^{[2]}$ is positive semidefinite, we conclude 
that $\omega$ is parallel and $\B^{[2]}\omega=0$. 

Hence, $\omega$ has constant 
positive norm. After multiplying $\omega$ 
by a positive constant, we may assume that 
$\|\omega\|=\sqrt{2}$. By Lemma \ref{proalg}, 
the formula
\[
\langle JX,Y\rangle=\omega(X,Y),
\qquad X,Y\in TM,
\]
defines a global orthogonal 
almost complex structure $J$ compatible 
with the given orientation. 
Since $\omega$ is parallel, $J$ is also 
parallel. Thus, $(M^4,g,J)$ is 
a K\"ahler manifold. 

We claim that 
\begin{equation}\label{Jx}
\a_f(JX,JY)=\a_f(X,Y)
\quad\text{for all }X,Y\in TM.
\end{equation}

Fix an arbitrary point $x\in M^4$. 
By Lemma \ref{proalg}, there exists a 
positively oriented orthonormal basis 
$\{e_i\}_{i=1}^4$ of $T_xM$ such that 
\[
Je_1=e_2,\qquad Je_3=e_4,\qquad
\omega_x=\omega_1\wedge\omega_2+
\omega_3\wedge\omega_4,
\]
where $\{\omega_i\}_{i=1}^4$ is the 
dual basis.

Since $\omega$ is parallel and belongs 
to the kernel of $\B^{[2]}$, the 
four-dimensional curvature identity yields 
\[
R_{1313}+R_{1414}+R_{2323}+R_{2424}
-2R_{1234}=0.
\]
Therefore, the equality conditions in 
part $(ii)$ of Proposition \ref{propu} 
apply to this orthonormal basis. 
In particular, the second fundamental 
form $\a_f$ satisfies \eqref{c1}--\eqref{c2}, 
with $\varepsilon=+1$. 
Since $Je_1=e_2$ and $Je_3=e_4$, 
these conditions imply that 
\[
\a_f(Je_i,Je_j)=\a_f(e_i,e_j)
\quad\text{ for all}\quad1\leq i,j\leq4.
\]
By bilinearity, \eqref{Jx} holds at $x$. 
Since $x$ was arbitrary, the immersion 
$f$ is $(2,0)$-geodesic with respect to 
the K\"ahler structure $J$. 

The case $b_-(M)>0$ is handled 
similarly. Indeed, a nonzero harmonic 
anti-self-dual $2$-form is parallel and, 
after normalization, defines an orthogonal 
parallel complex structure $J$ inducing 
the orientation opposite to the given one. 

By the anti-self-dual version of 
Lemma \ref{proalg}, at every point 
$x\in M^4$ there exists a positively 
oriented orthonormal basis 
$\{e_i\}_{i=1}^4$ such that 
$Je_1=e_2$ and $Je_3=-e_4$. 
The corresponding curvature identity is 
\[
R_{1313}+R_{1414}+R_{2323}+R_{2424}
+2R_{1234}=0.
\]
Applying the equality conditions of 
Proposition \ref{propu} with 
$\varepsilon=-1$, we again obtain 
\[
\a_f(JX,JY)=\a_f(X,Y)
\quad\text{ for all}\quad X,Y\in TM.
\]
Thus, $M^4$ admits a K\"ahler structure 
with respect to which $f$ is 
$(2,0)$-geodesic. 

Finally, suppose that the scalar 
curvature $\tau$ of $M^4$ vanishes 
identically. Since equality holds 
in \eqref{1} at every point, the 
Gauss equation gives 
\begin{align*}
0=\tau
=16H^2-S+
2\sum_{1\leq i<j\leq4}K_{\bar M}(f_*e_i\wedge f_*e_j)
\geq16H^2-S+12c=S+4c.
\end{align*}
This implies 
$S=c=0$. Consequently, $f$ is 
totally geodesic. Moreover, equality 
in the preceding estimate implies 
that all tangential sectional curvatures 
of $\bar M$ vanish along $f$. 
Thus, the Gauss equation shows that 
$M^4$ is flat. 

By the Gauss--Bonnet--Chern theorem, 
the Euler characteristic of $M^4$ 
vanishes. Since $M^4$ is compact, 
connected, and oriented, Poincar\'e 
duality gives 
\[
0=\chi(M)=2-2b_1(M)+b_2(M).
\]
Consequently, 
\[
b_2(M)=2b_1(M)-2.
\]
By hypothesis, $b_2(M)>0$. 
It follows that $b_1(M)\geq2$. 

On the other hand, since $M^4$ is flat, 
the Bochner formula implies that every 
harmonic $1$-form is parallel. Hence, 
$b_1(M)\leq4$. If $b_1(M)=4$, 
there exist four linearly independent 
parallel $1$-forms on $M^4$. Thus, 
the holonomy group 
is trivial, and $M^4$ is isometric 
to a flat torus. Conversely, every 
flat four-torus has first Betti 
number equal to $4$. 
Therefore, $2\leq b_1(M)\leq4$, and 
$b_1(M)=4$ if and only if $M^4$ 
is isometric to a flat torus. 
This completes the proof. 
\end{proof}

\subsection{Four-dimensional submanifolds with finite 
fundamental group}

\begin{theorem}\label{n=4p1f}
Let $f\colon M^4\to\bar {M}^{4+m}$ 
be an isometric immersion of a compact, 
oriented Riemannian four-manifold with finite 
fundamental group. Suppose that inequality 
\eqref{1} holds at every point, where $c\geq0$ 
is a lower bound on the curvature operator of 
$\bar {M}^{4+m}$. 
Then one of the following assertions holds: 
\begin{enumerate}[
label=(\roman*),
topsep=1pt,
itemsep=0.2ex,
leftmargin=*,
partopsep=0.2ex,
parsep=0.2ex,
]

\item The manifold $M^4$ is diffeomorphic 
to $\Sf^4$.

\item Equality holds in \eqref{1} at 
every point, and one of the following occurs: 
\begin{enumerate}[
label=(\alph*),
topsep=1pt,
itemsep=0.2ex,
partopsep=0.1ex,
parsep=0.2ex,
leftmargin=*
]

\item The manifold $M^4$ admits a K\"ahler 
structure with respect to which $f$ is 
$(2,0)$-geodesic. Moreover, $M^4$ is 
biholomorphic to $\CP^2$. 

\item The universal cover $\tilde M^4$ of 
$M^4$ is isometric to a Riemannian product 
$(N_1,g_1)\times(N_2,g_2)$, 
where each factor is diffeomorphic to 
$\Sf^2$ and has nonnegative Gaussian 
curvature. Moreover, there exist smooth 
Dupin principal normals 
$\delta_1$ and $\delta_2$ of $f$ 
satisfying 
\[
\pi_*(TN_i)\subset E_{\delta_i}, 
\quad i=1,2,\qquad
\langle\delta_1,\delta_2\rangle=-c,
\]
where $\pi\colon\tilde M^4\to M^4$ 
is the universal covering map. 
\end{enumerate}
\end{enumerate}
\end{theorem}

\begin{proof}
By Proposition \ref{propu}, the manifold 
$M^4$ has nonnegative isotropic curvature. 
We first establish the theorem when 
$M^4$ is simply connected. 
\vspace{1ex}

\noindent\emph{Case I.}
Suppose that $M^4$ has positive isotropic 
curvature at some point. 
By Remark $(iv)$ in \cite{ses}, 
the manifold $M^4$ admits a metric 
of positive isotropic curvature. 
Since $M^4$ is simply connected, 
the theorem of Micallef and Moore 
\cite{MM} implies that $M^4$ is 
homeomorphic to $\Sf^4$. 

We claim that $M^4$ is locally 
irreducible. Otherwise, since $M^4$ 
is compact and simply connected, 
Theorem \ref{duke} implies that 
$M^4$ is isometric to a Riemannian 
product of two compact surfaces. 
However, such a product admits, 
at every point, an orthonormal 
four-frame with vanishing isotropic 
curvature. This contradicts the 
assumption that $M^4$ has positive 
isotropic curvature at some point. 

Thus, $M^4$ is locally irreducible, 
and Theorem \ref{acta} applies. 
According to that theorem, one 
of the following occurs: 
\begin{enumerate}[
label=(\alph*),
topsep=1pt,
itemsep=0.2ex,
partopsep=0.2ex,
parsep=0.2ex,
leftmargin=*
]

\item $M^4$ is diffeomorphic to a 
spherical space form; 
\item $M^4$ is a K\"ahler manifold 
biholomorphic to $\CP^2$; 
\item $M^4$ is isometric to a 
compact symmetric space. 
\end{enumerate}

Since $M^4$ is homeomorphic to 
$\Sf^4$, the second alternative 
is excluded, because 
$b_2(\Sf^4)=0$, 
while $b_2(\CP^2)=1$. 
In the first alternative, since 
$M^4$ is simply connected, it is 
diffeomorphic to $\Sf^4$. 

In the third alternative, since 
$M^4$ is simply connected and 
locally irreducible, the 
classification of compact 
symmetric spaces in dimension 
four implies that $M^4$ is 
isometric, up to scaling, to 
either $\Sf^4$ or $\CP^2$. 
The latter possibility is excluded 
by the preceding Betti-number 
argument. Consequently, $M^4$ 
is diffeomorphic to $\Sf^4$. 
Thus, alternative $(i)$ holds. 
\vspace{1ex}

\noindent\emph{Case II.}
Suppose now that $M^4$ has no 
points of positive isotropic 
curvature. 
By part $(ii)$ of Proposition 
\ref{propu}, equality holds in 
\eqref{1} at every point. 
We distinguish two subcases. 
\vspace{0.5ex}

\noindent\emph{Subcase $II_a$.} 
Suppose that $M^4$ is locally 
irreducible. 
By Theorem \ref{acta}, one of 
the three alternatives described 
in \emph{Case I} occurs. 

If $M^4$ is diffeomorphic to a 
spherical space form, then, 
since it is simply connected, 
it is diffeomorphic to $\Sf^4$. 
Thus, alternative $(i)$ holds. 

Suppose next that $M^4$ is a 
K\"ahler manifold biholomorphic 
to $\CP^2$. 
Since $b_2(M)=1$, either 
$b_+(M)=1$ or $b_-(M)=1$. 
Therefore, Theorem \ref{n=4beta} 
implies that $M^4$ admits a 
K\"ahler structure with respect 
to which $f$ is $(2,0)$-geodesic. 
Hence, alternative $(ii)$-$(a)$ 
holds. 

Finally, suppose that $M^4$ is 
isometric to a compact symmetric 
space. Since $M^4$ is simply 
connected and locally irreducible, 
it is isometric, up to scaling, 
to either $\Sf^4$ or $\CP^2$. 
In the spherical case, $M^4$ 
is diffeomorphic to $\Sf^4$, 
and alternative $(i)$ holds. 
In the complex projective case, 
we have $b_2(M)=1$, and 
Theorem \ref{n=4beta} implies 
that $f$ is $(2,0)$-geodesic 
with respect to a K\"ahler 
structure on $M^4$. 
Thus, alternative $(ii)$-$(a)$ 
holds. 
\vspace{0.5ex}

\noindent\emph{Subcase $II_b$.} 
Suppose that $M^4$ is locally 
reducible. 
Since $M^4$ is compact and simply connected, 
Theorem \ref{duke} implies that $M^4$ is 
isometric to a Riemannian product 
$(N_1^2,g_1)\times(N_2^2,g_2)$ of compact, 
simply connected surfaces. 
By virtue of \eqref{1}, part $(iii)$ 
of Proposition \ref{propprod} implies that 
$M^4$ has nonnegative curvature operator. 
Consequently, both factors have nonnegative 
Gaussian curvature and are diffeomorphic 
to $\Sf^2$. 

Since $f$ satisfies \eqref{1}, 
part $(iii)$ of Proposition 
\ref{propprod} implies that 
there exist smooth 
Dupin principal normals 
$\delta_1$ and $\delta_2$ of $f$ 
satisfying 
\[
TN_i\subset E_{\delta_i}\quad
i=1,2,\qquad
\langle\delta_1,\delta_2\rangle=-c.
\]
Consequently, the product has the 
geometry described in 
alternative $(ii)$-$(b)$. 
This completes the proof 
when $M^4$ is simply connected. 
\vspace{1ex}

We now consider the case where the 
fundamental group $\pi_1(M^4)$ 
is finite. Let $\pi\colon\tilde M^4\to M^4$ 
be the universal covering map, and set 
$\tilde f=f\circ\pi$. Since $M^4$ is compact 
and $\pi_1(M^4)$ is finite, the universal 
cover $\tilde M^4$ is compact. Moreover, 
$\tilde f$ satisfies inequality \eqref{1} 
at every point. By the preceding argument, 
one of the following occurs: 
\begin{enumerate}[
label=(\alph*),
topsep=1pt,
itemsep=0.2ex,
partopsep=0.2ex,
parsep=0.2ex,
leftmargin=*
]
\item $\tilde M^4$ is diffeomorphic 
to $\Sf^4$; 
\item $\tilde M^4$ is biholomorphic 
to $\CP^2$ and admits a K\"ahler 
structure with respect to which 
$\tilde f$ is $(2,0)$-geodesic; 
\item $\tilde M^4$ is isometric 
to a Riemannian product 
$(N_1,g_1)\times(N_2,g_2)$, 
where each factor is diffeomorphic 
to $\Sf^2$ and has nonnegative 
Gaussian curvature, and $\tilde f$ 
has the properties established 
in Subcase $II_b$.
\end{enumerate}

Suppose first that either (a) or (b) 
holds. Let $d=|\pi_1(M^4)|$ denote the 
number of sheets of the universal covering. 
Since $\pi_1(M^4)$ is finite, we have 
$b_1(M)=0$. By Poincar\'e 
duality, 
\[
\chi(M)=2+b_2(M)\geq2.
\]
On the other hand, the multiplicativity 
of the Euler characteristic under finite 
coverings gives 
\[
\chi(\tilde M)=d\,\chi(M).
\]
Since $\chi(\tilde M)$ is either 
$2$ or $3$, it follows that $d=1$. 
Thus, $\pi$ is a diffeomorphism, 
and either alternative $(i)$ or 
alternative $(ii)$-$(a)$ holds. 

Finally, suppose that $(c)$ holds. 
Then $\tilde M^4$ is isometric to 
$(N_1,g_1)\times(N_2,g_2)$, where 
each factor is diffeomorphic to 
$\Sf^2$ and has nonnegative Gaussian 
curvature. Moreover, there exist 
smooth Dupin principal normals 
$\tilde\delta_1$ and $\tilde\delta_2$ 
of $\tilde f$ satisfying 
\[
TN_i\subset E_{\tilde\delta_i},
\quad i=1,2,
\qquad
\langle\tilde\delta_1,\tilde\delta_2\rangle=-c.
\]

We claim that every deck transformation 
preserves each of the two product factors. 
Indeed, by the uniqueness of the de Rham 
decomposition, every deck transformation 
either preserves the factors individually 
or interchanges them. 

Suppose that a nontrivial deck transformation 
$\gamma$ interchanges the factors. 
By the K\"unneth formula, 
\[
H^2(\tilde M;\R)
=p_1^*H^2(N_1;\R)
\oplus p_2^*H^2(N_2;\R),
\]
where $p_i$ denotes the projection 
onto $N_i$. Since $\gamma$ interchanges 
the factors, its induced action on 
$H^2(\tilde M;\R)$ interchanges the 
two one-dimensional summands. Hence, 
\[
\operatorname{tr}
\bigl(\gamma^*|_{H^2(\tilde M;\R)}\bigr)=0.
\]
Since $M^4$ is oriented, every deck 
transformation preserves orientation. 
Moreover, 
$H^1(\tilde M;\R)=
H^3(\tilde M;\R)=0$. Consequently, 
the Lefschetz number of $\gamma$ is 
$L(\gamma)=2$. 
The Lefschetz fixed-point theorem then 
implies that $\gamma$ has a fixed point, 
contradicting the fact that nontrivial 
deck transformations act freely. 
Hence, every deck transformation 
preserves the tangent distributions 
$TN_1$ and $TN_2$. Since the second 
fundamental form of $\tilde f$ is the 
pullback of that of $f$, the principal 
normals $\tilde\delta_1$ and 
$\tilde\delta_2$ are invariant under 
the deck transformations. Therefore, 
they descend to smooth Dupin 
principal normals $\delta_1$ and 
$\delta_2$ of $f$ satisfying 
\[
\pi_*(TN_i)\subset E_{\delta_i},
\quad i=1,2,
\qquad
\langle\delta_1,\delta_2\rangle=-c.
\]
Thus, alternative $(ii)$-$(b)$ holds, 
and the proof is complete. 
\end{proof}

\begin{remark}
\emph{In alternative $(ii)$-$(b)$ of 
Theorem \ref{n=4p1f}, the 
fundamental group of $M^4$ 
is either trivial or isomorphic 
to $\mathbb Z_2$. 
Indeed, let $d=|\pi_1(M^4)|$ be the 
number of sheets of the universal covering. 
Since the universal cover is 
diffeomorphic to $\Sf^2\times\Sf^2$, 
we have $\chi(\tilde M)=4$. 
As $M^4$ is compact, 
oriented, and has finite 
fundamental group, Poincar\'e 
duality gives 
\[
\chi(M)=2+b_2(M)\geq2.
\]
Consequently, 
\[
4=d(2+b_2(M))\geq2d,
\]
and hence $d\leq2$. Thus, 
$\pi_1(M^4)=\{1\}$ or 
$\pi_1(M^4)\cong\Z_2$. 
In particular, if $\pi_1(M^4)$ 
is nontrivial, then $d=2$, 
$b_2(M)=0$, and 
$\chi(M)=2$. 
}
\end{remark}

\subsection{Four-dimensional submanifolds with infinite fundamental group}
We next establish the classification of 
four-dimensional submanifolds with infinite 
fundamental group satisfying \eqref{1}. 

\begin{theorem}\label{n=4p1inf}
Let $f\colon M^4\to\bar {M}^{4+m}$ 
be an isometric immersion of a compact, 
oriented Riemannian four-manifold 
with infinite fundamental group. Suppose that 
inequality \eqref{1} holds at every point, 
where $c\geq0$ is a lower bound on the 
curvature operator of $\bar {M}^{4+m}$. 
Then one of the following assertions holds: 
\begin{enumerate}[
label=(\roman*),
topsep=1pt,
itemsep=0.2ex,
leftmargin=*,
partopsep=0.2ex,
parsep=0.2ex,
]

\item The universal cover of $M^4$ 
is isometric to a Riemannian product 
$\R\times N^3$, where $N^3$ is 
diffeomorphic to $\Sf^3$ and has 
nonnegative Ricci curvature. 

\item Equality holds in \eqref{1} 
at every point, $c=0$, and one of the 
following occurs: 
\begin{enumerate}[
label=(\alph*),
topsep=1pt,
itemsep=0.2ex,
partopsep=0.2ex,
parsep=0.2ex,
leftmargin=*
]
\item The universal cover $\tilde M^4$ 
of $M^4$ is isometric to a Riemannian 
product $(\R^2,g_0)\times(N^2,g_N)$, 
where $g_0$ is the Euclidean metric 
and $N^2$ is diffeomorphic to $\Sf^2$ 
with nonnegative Gaussian curvature. 
Moreover, $\pi_*(T\R^2)\subset\mathcal D_f$, 
and the normal vector field $\delta=2\mathcal H_f$ 
is a Dupin principal normal satisfying 
$\pi_*(TN)\subset E_\delta$, 
where $\pi\colon\tilde M^4\to M^4$ 
is the universal covering map. 

\item The manifold $M^4$ is flat, 
with $1\leq b_1(M)\leq4$, 
and the immersion $f$ is totally geodesic. 
\end{enumerate}
\end{enumerate}
\end{theorem}
\proof
Proposition \ref{propu} implies that $M^4$ has 
nonnegative isotropic curvature. We first claim 
that $M^4$ is locally reducible. Suppose, 
to the contrary, that $M^4$ is locally irreducible. 
Then one of the cases $(i)$--$(iii)$ in Theorem 
\ref{acta} applies. In each of these cases, 
the universal cover $\tilde M^4$ is compact, 
which contradicts the assumption that the 
fundamental group of $M^4$ is infinite. 

Hence, $M^4$ is locally reducible. It then follows 
from Theorem \ref{duke} that the universal cover 
$\tilde M^4$ is isometric to a Riemannian product 
as in cases $(a)$ or $(b)$ of that theorem. 
The proof of the theorem is divided into two cases.
\vspace{0.5ex}

\indent\emph{Case I}. Suppose that $M^4$ has 
positive isotropic curvature at some point. Since 
$M^4$ has nonnegative isotropic curvature 
everywhere, Remark (iv) in \cite{ses} implies 
that $M^4$ admits a metric of positive 
isotropic curvature. By the theorem of Micallef and 
Moore \cite{MM}, we have $\pi_2(M^4)=0$. 
Since covering maps induce isomorphisms on 
all higher homotopy groups (see, for instance, 
Proposition 4.1 in \cite{Hat}), it follows that 
$\pi_2(\tilde M^4)=0$. On the other hand, 
$\pi_2(\Sf^2)\cong\Z$, 
which rules out case~(b) of 
Theorem~\ref{duke}, since the corresponding 
product has a factor diffeomorphic to 
$\Sf^2$. 

Moreover, the flat case is excluded, 
since $M^4$ has positive isotropic 
curvature at some point. 
Hence, the universal cover $\tilde M^4$ is 
isometric to a Riemannian product 
$(\R^{n_0},g_0)\times(N_1^{n_1},g_1)$, 
where $n_0\geq1$, $g_0$ is the Euclidean 
metric, and either 
$n_1=2$ and $N_1$ is diffeomorphic to $\Sf^2$ 
with nonnegative Gaussian curvature, or 
$n_1=3$ and $N_1$ is compact with 
nonnegative Ricci curvature. 
Since $\pi_2(\tilde M^4)=0$, the former
possibility is also excluded. Thus,
$n_0=1$ and $n_1=3$, so that
$\tilde M^4$ is isometric to a
Riemannian product $\R\times N_1$,
where $N_1$ is a compact, simply connected
Riemannian three-manifold with nonnegative
Ricci curvature. By Theorem~1.1
of \cite{H3}, $N_1$ is diffeomorphic
to $\Sf^3$.
\vspace{0.5ex}

\indent\emph{Case II}. Now suppose that 
$M^4$ has no point of positive 
isotropic curvature. Then part~(ii) of 
Proposition~\ref{propu} implies that equality 
holds in \eqref{1} at every point. 

We first claim that case~(b) of 
Theorem~\ref{duke} cannot occur. 
Suppose, to the contrary, that $\tilde M^4$ 
is isometric to a Riemannian product 
$(N_1,g_1)\times(N_2,g_2)$, where $N_2$ 
is compact and the Gaussian curvatures 
of the two factors satisfy 
$\min K_{N_2}\geq-\min K_{N_1}>0$. 
Clearly, the isometric immersion $\tilde f=f\circ\pi$ 
satisfies condition \eqref{1}, where 
$\pi\colon\tilde M^4\to M^4$ is the 
covering map. By part~$(iii)$ of 
Proposition~\ref{propprod}, both factors 
have nonnegative Gaussian curvature, 
a contradiction.

Consequently, only case $(a)$ of 
Theorem \ref{duke} can occur. 
Thus, the universal cover $\tilde M^4$ 
is isometric to one of the products 
$\R\times N^3$, $\R^2\times N^2$, 
$\R^4$, where $N^3$ is compact, simply 
connected, and has nonnegative Ricci 
curvature, while $N^2$ is diffeomorphic 
to $\Sf^2$ and has nonnegative 
Gaussian curvature. 

If $\tilde M^4$ is isometric to 
$\R\times N^3$, then Theorem~1.1 
of \cite{H3} implies that $N^3$ 
is diffeomorphic to $\Sf^3$. 
Hence, alternative~$(i)$ holds. 

Suppose next that $\tilde M^4$ is 
isometric to $\R^2\times N^2$. 
By part $(iii)$ of Proposition 
\ref{propprod}, the second fundamental 
form of $\tilde f=f\circ\pi$ is adapted 
to the product net. Moreover, there 
exist smooth Dupin principal normals 
$\tilde\delta_1$ and $\tilde\delta_2$ 
satisfying 
\[
T\R^2\subset E_{\tilde\delta_1},
\qquad
TN^2\subset E_{\tilde\delta_2},
\quad\text{and}\quad
\langle\tilde\delta_1,\tilde\delta_2\rangle=-c.
\]

Let $v_1,v_2\in T\R^2$ be orthonormal. 
Since the plane spanned by $v_1$ and 
$v_2$ has zero sectional curvature, 
the Gauss equation yields 
\[
0=K_{\bar M}
(\tilde f_*v_1\wedge\tilde f_*v_2)
+\|\tilde\delta_1\|^2
\geq c+\|\tilde\delta_1\|^2.
\]
Hence, $c=0$ and $\tilde\delta_1=0$. 
Consequently, 
$T\R^2\subset\mathcal D_{\tilde f}$. 
Moreover, since the second fundamental 
form is adapted to the product net, 
we have 
$\tilde\delta_2=2\mathcal H_{\tilde f}$. 
Since $\pi$ is a local isometry, 
these conclusions descend to $M^4$. 
In particular, the normal vector field 
$\delta=2\mathcal H_f$ 
is a Dupin principal normal satisfying 
$\pi_*(T\R^2)\subset\mathcal D_f$, and 
$\pi_*(TN^2)\subset E_\delta$. 
Therefore, alternative $(ii)$-$(a)$ 
holds. 

Finally, suppose that $\tilde M^4=\R^4$. 
Then $M^4$ is a compact, oriented flat 
manifold. By the Gauss--Bonnet--Chern 
theorem, $\chi(M)=0$. Since $M^4$ is 
compact and oriented, Poincar\'e duality 
gives 
\[
0=\chi(M)=2-2b_1(M)+b_2(M).
\]
Consequently, 
\[
b_2(M)=2b_1(M)-2,
\]
and hence $b_1(M)\geq1$. 
On the other hand, every harmonic 
$1$-form on a compact flat manifold 
is parallel. Therefore, 
$b_1(M)\leq4$, with equality 
if and only if $M^4$ is isometric 
to a flat torus. 

Since the scalar curvature 
$\tau$ of $M^4$ vanishes identically 
and equality holds in \eqref{1}, given an 
orthonormal basis of the tangent space at 
any point, the Gauss equation yields 
\begin{align*}
0=\tau=16H^2-S+
2\sum_{1\leq i<j\leq4}
K_{\bar M}(f_*e_i\wedge f_*e_j)
\geq16H^2-S+12c=S+4c.
\end{align*}
Consequently, $S=c=0$. 
Thus, the immersion $f$ is 
totally geodesic, and alternative 
$(ii)$-$(b)$ holds. 
This completes the proof. 
\qed
\vspace{1ex}

\noindent\emph{Proof of Theorem \ref{n=4}:} 
The assertion follows directly from Theorems 
\ref{n=4p1f} and \ref{n=4p1inf}.
\qed

\section{Proof of the main result for the case $n\geq5$}
The proof of Theorem \ref{n5} is divided into 
two cases, depending on whether the 
fundamental group of the submanifold 
is finite or infinite. 

Before proceeding, we first establish the 
following result. 

\begin{theorem}\label{ceven}
Let $f\colon M^n\to\bar M^{n+m}$, 
$n=2k\geq6$, be an isometric immersion 
of a compact, locally irreducible 
Riemannian manifold satisfying inequality 
\eqref{1} at every point, where $c$ is a 
lower bound on the curvature operator of 
$\bar M^{n+m}$. 
If $b_2(M^n)>0$, then $c\geq0$, 
the immersion $f$ is totally geodesic, 
and $M^n$ is either a compact, simply 
connected, irreducible Hermitian symmetric 
space or a K\"ahler manifold biholomorphic 
to $\CP^k$. 
\end{theorem}
\proof
By Proposition \ref{propu}, the manifold $M^n$ 
has nonnegative isotropic curvature. Moreover, 
part $(i)$ of Proposition \ref{propeven} implies 
that the Bochner--Weitzenb\"ock operator 
$\B^{[2]}$ is positive semidefinite. Since 
$b_2(M^n)>0$, there exists a 
nonzero harmonic 2-form 
$\omega\in\Omega^2(M^n)$. By the 
Bochner--Weitzenb\"ock formula, $\omega$ 
is parallel and $\<\B^{[2]}\omega,\omega\>=0$ 
everywhere. In particular, $\omega$ is 
nowhere vanishing, and part $(ii)$ of 
Proposition \ref{propeven} applies at 
every point. 

Define the skew-symmetric endomorphism 
$L$ of $TM$ by 
\[
\omega(X,Y)=\<LX,Y\>,
\qquad X,Y\in TM.
\]
Since $\omega$ is parallel, so are $L$ 
and the self-adjoint endomorphism $L^2$. 
Consequently, the eigenspaces of $L^2$ 
corresponding to distinct eigenvalues 
would determine an orthogonal splitting 
of $TM$ into parallel subbundles. 
Since $M^n$ is locally irreducible, 
we conclude that $L^2=\mu\,{\rm Id}$ 
for some constant $\mu\in\R$. 

At any point $x\in M^n$, choose an 
orthonormal basis $\{e_i\}_{i=1}^n$ 
of $T_xM$ such that $e_{2a-1}=v_a$ 
and $e_{2a}=w_a$ for $1\leq a\leq k$, 
where $\{v_1,w_1,\dots,v_k,w_k\}$ 
is chosen so that \eqref{Lambda} holds. 
Then 
\[
Le_{2a-1}=\lambda_a e_{2a},
\qquad
Le_{2a}=-\lambda_a e_{2a-1},
\quad 1\leq a\leq k.
\]
Since $\omega_x\neq0$ and 
$L^2=\mu\,{\rm Id}$, we obtain 
\be\label{mu<0}
\lambda_1^2=\cdots=\lambda_k^2=-\mu>0.
\ee
Hence, 
$
J=\frac{1}{\sqrt{-\mu}}L
$
defines an orthogonal almost complex 
structure on $M^n$. Moreover, since 
$L$ is parallel, so is $J$. Thus, 
$(M^n,\<\cdot,\cdot\>,J)$ is a 
K\"ahler manifold. It follows from 
Theorem 2.1 of \cite{MW} that $M^n$ 
is simply connected. 

We next claim that $f$ is totally 
geodesic. Otherwise, there would exist 
a point $x_0\in M^n$ at which $f$ is 
not totally geodesic. By part $(ii)$-$(b)$ 
of Proposition \ref{propeven}, exactly 
one of the numbers $\lambda_a$, 
$1\leq a\leq k$, would be nonzero 
at $x_0$, contradicting \eqref{mu<0}. 
Therefore, $f$ is totally geodesic, 
and inequality \eqref{1} implies 
that $c\geq0$. 

Since $f$ is totally geodesic, the 
Gauss equation implies that $M^n$ 
has nonnegative curvature operator. 
Moreover, $M^n$ cannot be homeomorphic 
to $\Sf^n$, since $b_2(M^n)>0$. 
Thus, by Theorem 5 of \cite{chow}, 
$M^n$ is either a compact symmetric 
space or a K\"ahler manifold 
biholomorphic to $\CP^k$. 
In the former case, since $M^n$ 
is K\"ahler and locally irreducible, 
it is a compact simply connected 
irreducible Hermitian symmetric space.\qed

\subsection{Submanifolds with finite fundamental group}

The aim of this section is to establish the 
following result for submanifolds with 
finite fundamental group satisfying the 
pinching condition \eqref{1}. 

\begin{theorem}\label{p1fn5}
Let $f\colon M^n\to\bar {M}^{n+m}$, 
$n\geq5$, be an isometric immersion of a compact 
Riemannian manifold with finite fundamental group. 
Suppose that inequality \eqref{1} holds at every 
point, where $c\geq0$ is a lower bound on the 
curvature operator of $\bar {M}^{n+m}$. 
Then one of the following alternatives holds: 
\begin{enumerate}[
label=(\roman*),
topsep=1pt,
itemsep=0.2ex,
partopsep=0.2ex,
parsep=0.2ex,
leftmargin=*
]

\item The manifold $M^n$ is 
diffeomorphic to a spherical space form. 
\item Equality holds in \eqref{1} at every 
point, and $c=0$. Moreover, one of the 
following occurs: 
\begin{enumerate}[
label=(\alph*),
topsep=1pt,
itemsep=0.2ex,
partopsep=0.2ex,
parsep=0.2ex,
leftmargin=*
]
\item The immersion $f$ is totally geodesic, 
and the universal cover of $M^n$ is either 
a simply connected irreducible compact 
symmetric space or a K\"ahler manifold 
biholomorphic to $\CP^{n/2}$. 
\item The immersion $f$ is totally geodesic, 
and the universal cover of $M^n$ is isometric 
to a Riemannian product 
\[
(N_1^{n_1},g_1)\times\cdots\times
(N_r^{n_r},g_r),\qquad r\geq2,
\]
where each factor $(N_i,g_i)$ is 
irreducible and satisfies 
$2\leq n_i<n-2$ for all $1\leq i\leq r$. 
\item The universal cover $\tilde M^n$ of 
$M^n$ is isometric to a Riemannian product 
$(N_1^2,g_1)\times(N_2^{n-2},g_2)$, 
where $(N_2^{n-2},g_2)$ is irreducible. 
Moreover, $\pi_*(TN_1)\subset\mathcal D_f$, 
and the normal vector field 
$\delta=n\mathcal H_f/(n-2)$ 
is a Dupin principal normal of $f$ 
satisfying $\pi_*(TN_2)\subset E_\delta$, 
where $\pi\colon\tilde M^n\to M^n$ 
is the universal covering map. 
\end{enumerate}
In each of the cases $(a)$--$(c)$, 
every irreducible factor of the 
universal cover having dimension at 
least three is either a compact, 
simply connected, irreducible 
symmetric space, a Riemannian 
manifold with positive 
curvature operator diffeomorphic to 
a sphere, or a K\"ahler manifold 
biholomorphic to a complex 
projective space. Every two-dimensional 
factor is diffeomorphic to $\Sf^2$ with 
nonnegative Gaussian curvature. 
\end{enumerate}
\end{theorem}
\proof
By Proposition \ref{propu}, the manifold $M^n$ 
has nonnegative isotropic curvature. We first 
prove the theorem under the additional assumption 
that $M^n$ is simply connected. We distinguish 
two cases. 
\vspace{1ex}

\noindent\emph{Case I.}
Suppose that $M^n$ has positive isotropic 
curvature at some point. Since $M^n$ has 
nonnegative isotropic curvature, Remark $(iv)$ 
in \cite{ses} implies that $M^n$ admits a 
metric of positive isotropic curvature. 
By the theorem of Micallef and Moore 
\cite{MM}, $M^n$ is homeomorphic to $\Sf^n$. 

We claim that $M^n$ is locally irreducible. 
Otherwise, since $M^n$ is complete and simply 
connected, the de Rham decomposition theorem 
would yield a nontrivial isometric splitting 
$M^n=M_1\times M_2$, where both factors 
have positive dimension. This is impossible, 
since $\Sf^n$ cannot be expressed as a 
topological product of two positive-dimensional 
compact manifolds. Hence, $M^n$ is locally 
irreducible. 

Therefore, $M^n$ falls under the alternatives 
described in Theorem \ref{acta}. Since $M^n$ 
is homeomorphic to $\Sf^n$, case $(ii)$ 
of that theorem is excluded. In case $(i)$, 
$M^n$ is diffeomorphic to $\Sf^n$. In the remaining 
case $(iii)$, $M^n$ is isometric to a compact 
irreducible symmetric space. Hence $M^n$ is 
isometric to a round sphere. 
This proves alternative $(i)$. 
\vspace{1ex}

\noindent\emph{Case II.} 
Suppose now that $M^n$ has no points of 
positive isotropic curvature. It follows 
from part $(ii)$ of Proposition \ref{propu} 
that equality holds in \eqref{1} at every 
point. We now distinguish two subcases. 
\vspace{0.5ex}

\noindent\emph{Subase $II_a$.} 
We first assume that $M^n$ is locally 
reducible. By the de Rham decomposition 
theorem, $M^n$ is isometric to a nontrivial 
Riemannian product. Moreover, by Theorem 
\ref{duke}, its decomposition must be 
of one of the two types described in 
that theorem. 

Case $(b)$ of Theorem \ref{duke} is ruled out 
by Proposition \ref{propprod}, which implies 
that $M^n$ has nonnegative curvature 
operator. Thus, $M^n$ is a Riemannian 
product as in case $(a)$ of Theorem 
\ref{duke}. Since each factor 
in that decomposition has dimension at 
least two, Proposition \ref{propprod} applies. 
For dimension reasons, only cases $(i)$ and 
$(ii)$ of that proposition can occur. In either 
case, Proposition \ref{propprod} implies that 
$c=0$ and that $M^n$ has nonnegative 
curvature operator. 

Suppose first that case $(i)$ of 
Proposition \ref{propprod} occurs. 
Then $f$ is totally geodesic. Moreover, 
$M^n$ is isometric to a Riemannian 
product 
\[
(N_1^{n_1},g_1)\times\cdots\times
(N_r^{n_r},g_r),\qquad r\geq2,
\]
where each factor $(N_i,g_i)$ is compact, 
simply connected, irreducible and satisfies 
$2\leq n_i<n-2$ for all $1\leq i\leq r$. 
By Theorem 5 of \cite{chow}, each 
irreducible factor is one of the following: 
\begin{enumerate}[
label=(\alph*),
topsep=1pt,
itemsep=0.2ex,
partopsep=0.2ex,
parsep=0.2ex,
leftmargin=*
]
 
\item[(1)] a simply connected compact 
symmetric space; 
\item[(2)] a Riemannian manifold 
admitting a metric with positive curvature 
operator; 
\item[(3)] a K\"ahler manifold 
biholomorphic to a complex 
projective space. 
\end{enumerate}
In case $(2)$, the factor is 
diffeomorphic to a sphere by 
the theorem of B\"ohm and Wilking 
\cite{bohm}. Hence, alternative 
$(ii)$-$(b)$ holds. 

Suppose next that case $(ii)$ of 
Proposition \ref{propprod} occurs. 
Then $M^n$ is isometric to a 
Riemannian product 
$(N_1^2,g_1)\times(N_2^{n-2},g_2)$, 
where $(N_1^2,g_1)$ is diffeomorphic to 
$\Sf^2$ with nonnegative Gaussian curvature 
and $(N_2^{n-2},g_2)$ is irreducible. 

As above, Theorem 5 of \cite{chow}, 
together with the theorem of B\"ohm 
and Wilking \cite{bohm}, implies 
that $(N_2^{n-2},g_2)$ is either 
a simply connected compact symmetric 
space, a Riemannian manifold with 
positive curvature operator 
diffeomorphic to a sphere, or a 
K\"ahler manifold biholomorphic 
to a complex projective space. 
Furthermore, by Proposition 
\ref{propprod}, we have 
$TN_1\subset\mathcal D_f$, 
and the normal vector field 
$\delta=n\mathcal H_f/(n-2)$ 
is a Dupin principal normal of $f$ 
satisfying $TN_2\subset E_\delta$. 
Thus, alternative $(ii)$-$(c)$ holds. 
%\vspace{0.1ex} 
\vspace{-10pt}

\noindent\emph{Subase $II_b$.} 
We now assume that $M^n$ is locally 
irreducible. Then Theorem \ref{acta} 
applies. 

In case $(i)$ of that theorem, 
$M^n$ is diffeomorphic to $\Sf^n$, 
and alternative $(i)$ follows. 

Suppose that case $(ii)$ occurs. 
Then $n$ is even, and $M^n$ is a 
K\"ahler manifold biholomorphic 
to $\CP^{n/2}$. In particular, 
$b_2(M^n)=1$. 
Since $M^n$ is compact and locally 
irreducible, Theorem \ref{ceven} 
implies that $f$ is totally geodesic. 
Moreover, since equality holds 
in \eqref{1} at every point, we 
conclude that $c=0$. Hence, 
alternative $(ii)$-$(a)$ holds. 

Finally, suppose that case $(iii)$ 
of Theorem \ref{acta} occurs, so 
that $M^n$ is isometric to a compact 
irreducible symmetric space. 
Since $M^n$ has no points of positive 
isotropic curvature, its sectional 
curvature cannot be constant unless 
$M^n$ is flat. The latter possibility 
is excluded because $M^n$ is compact 
and simply connected. 
We claim that $f$ is totally geodesic. 
Indeed, Proposition~\ref{propu}(ii) and 
Proposition~\ref{fillup} together imply 
that $E_\delta(x)=T_xM$ at every point $x$,
where $\delta=n\mathcal H_f(x)/(n-2)$.
Thus, at every point $x\in M^n$, we have
\[
\alpha_f(X,Y)=\langle X,Y\rangle\delta,
\qquad X,Y\in T_xM.
\]
Taking the trace yields
$\mathcal H_f(x)=\delta$.
Since $\delta=n\mathcal H_f(x)/(n-2)$,
we conclude that $\delta=0$.
Hence, $f$ is totally geodesic.
Since equality holds in \eqref{1} 
at every point, it follows that 
$c=0$. Thus, alternative 
$(ii)$-$(a)$ holds. 

This completes the proof when 
$M^n$ is simply connected. 
\vspace{1ex}

We now consider the general case 
where $\pi_1(M^n)$ is finite. 
Let $\pi\colon\tilde M^n\to M^n$ 
be the universal covering map. 
Since $M^n$ is compact and its 
fundamental group is finite, 
$\tilde M^n$ is also compact. 
The composition 
$\tilde f=f\circ\pi\colon
\tilde M^n\to\bar {M}^{n+m}$ 
is an isometric immersion satisfying 
inequality \eqref{1} at every point. 
Moreover, $\tilde M^n$ is simply connected. 
We distinguish two cases. 
\vspace{1ex}

\noindent\emph{Case III.} 
Suppose that $M^n$ has positive isotropic 
curvature at some point. Since $M^n$ has 
nonnegative isotropic curvature, by Seshadri's 
\cite{ses} deformation result, $M^n$ admits a 
metric of positive isotropic curvature. The 
pullback of this metric to $\tilde M^n$ also 
has positive isotropic curvature. 
As in \emph{Case I}, we conclude that 
$\tilde M^n$ is diffeomorphic to $\Sf^n$. 
Hence, $M^n$ is locally 
irreducible and therefore one of the 
alternatives in Theorem \ref{acta} holds. 
In case $(i)$ of that theorem, $M^n$ is 
diffeomorphic to a spherical space form. 
Case $(ii)$ of that theorem is excluded since 
$\tilde M^n$ is diffeomorphic to $\Sf^n$. In 
the remaining case $(iii)$, $\tilde M^n$ is isometric
to a round sphere, giving alternative $(i)$. 
\vspace{1ex}

\noindent\emph{Case IV.} 
Suppose now that $M^n$ has no points of 
positive isotropic curvature. By part $(ii)$ 
of Proposition \ref{propu}, equality holds 
in \eqref{1} for $f$ at every point. Since $\pi$ is 
a local isometry, $\tilde M^n$ has no points 
of positive isotropic curvature, and $\tilde f$ 
satisfies equality in \eqref{1} 
at every point. 

We first suppose that $M^n$ is 
locally irreducible. Then Theorem \ref{acta} 
applies. In case $(i)$ of that theorem 
$M^n$ is diffeomorphic to a spherical 
space form. Suppose that either case 
$(ii)$ or case $(iii)$ holds. Clearly, $\tilde M^n$ is 
also locally irreducible. By the classification 
established in Subcase $IIb$, $\tilde f$ is 
totally geodesic, $c=0$, and $\tilde M^n$ 
is either a simply connected irreducible compact 
symmetric space or a K\"ahler manifold 
biholomorphic to $\CP^{n/2}$. 
Since $\pi$ is a local isometry and 
$\tilde f=f\circ\pi$, the immersion $f$ 
is also totally geodesic. Thus, 
alternative $(ii)$-$(a)$ holds. 

Suppose now that $M^n$ is 
locally reducible. Clearly, $\tilde M^n$ is 
also locally reducible. Then the 
preceding classification yields 
one of the following possibilities. 

First, $\tilde f$ is totally geodesic, 
and $\tilde M^n$ is isometric to a 
Riemannian product 
\[
(N_1^{n_1},g_1)\times\cdots\times 
(N_r^{n_r},g_r),\qquad r\geq2, 
\]
where each factor $(N_i,g_i)$ 
is compact, simply connected, irreducible 
and satisfies $2\leq n_i<n-2$, for all 
$1\leq i\leq r$. Clearly, $f$ 
is totally geodesic as well. 
Hence, alternative $(ii)$-$(b)$ 
holds. 

Alternatively, $\tilde M^n$ is 
isometric to a Riemannian product 
$(N_1^2,g_1)\times(N_2^{n-2},g_2)$, 
where $(N_1^2,g_1)$ is diffeomorphic 
to $\Sf^2$ with nonnegative Gaussian 
curvature, and $(N_2^{n-2},g_2)$ is irreducible. 
Furthermore, the second fundamental 
form of $\tilde f$ satisfies 
$TN_1\subset\mathcal D_{\tilde f}$, 
and 
$TN_2\subset E_{\tilde\delta}$, 
where 
$\tilde\delta=n\mathcal H_{\tilde f}/(n-2)$ 
is a Dupin principal normal of 
$\tilde f$. Since $\pi$ is a local isometry, 
the second fundamental forms of 
$f$ and $\tilde f$ are related by 
$$
\alpha_{\tilde f}(X,Y) 
=\alpha_f(\pi_*X,\pi_*Y),\qquad X,Y\in T\tilde M. 
$$
Consequently, 
$\pi_*(TN_1)\subset\mathcal D_f$, 
and $\pi_*(TN_2)\subset E_\delta$, 
where $\delta=n\mathcal H_f/(n-2)$ 
is a Dupin principal normal of $f$. 
Thus, alternative $(ii)$-$(c)$ holds. 

In both reducible cases, the 
classification of the compact irreducible 
factors follows from Theorem 5 
of \cite{chow} and the theorem 
of B\"ohm and Wilking \cite{bohm}. 
This completes the proof. 
\qed

\subsection{Submanifolds with infinite fundamental group}

We now consider compact submanifolds with infinite 
fundamental group satisfying the pinching condition 
\eqref{1}. More precisely, we prove the following 
result. 

\begin{theorem}\label{p1ifn5}
Let $f\colon M^n\to\bar {M}^{n+m}$, 
$n\geq5$, be an isometric immersion of a compact 
Riemannian manifold with infinite fundamental group. 
Suppose that inequality \eqref{1} holds at every point, 
where $c\geq0$ is a lower bound on the 
curvature operator of $\bar {M}^{n+m}$. Then 
one of the following occurs: 

\begin{enumerate}[
label=(\roman*),
topsep=1pt,
itemsep=0.2ex,
partopsep=0.2ex,
parsep=0.2ex,
leftmargin=*
]

\item The universal cover of $M^n$ is 
isometric to a Riemannian product 
$\R\times N^{n-1}$, 
where $N^{n-1}$ is diffeomorphic to $\Sf^{n-1}$ 
with nonnegative isotropic curvature. 

\item Equality holds in \eqref{1} at every 
point, $c=0$, and $f$ is totally geodesic. 
In particular, $M^n$ has nonnegative curvature 
operator, and its universal cover is isometric 
to a Riemannian product 
\[
(\R^{n_0},g_0)\times
(N_1^{n_1},g_1)\times\cdots\times
(N_r^{n_r},g_r),
\]
where $n_0\geq1$, $g_0$ is the Euclidean 
metric, $r\geq0$, and the following 
dimensional restrictions hold: 
\begin{enumerate}[
label=(\alph*),
topsep=1pt,
itemsep=0.2ex,
partopsep=0.2ex,
parsep=0.2ex,
leftmargin=*
]

\item If $r=1$, then
$1\leq n_0\leq n-2$ and $n_1=n-n_0$.
%If $r=1$, then either $(n_0,n_1)=(1,n-1)$ or $(n_0,n_1)=(n-2,2)$. 
\item If $r\geq2$, then 
$2\leq n_i<n-2$ for every 
$1\leq i\leq r$. 
\end{enumerate}
Each factor $(N_i,g_i)$ is compact, 
simply connected, irreducible, and 
has nonnegative curvature operator. 

\item Equality holds in \eqref{1} 
at every point, $c=0$, and the universal 
cover $\tilde M^n$ of $M^n$ is isometric 
to a Riemannian product 
$(\R^2,g_0)\times(N^{n-2},g_N)$, 
where $(N^{n-2},g_N)$ is compact, simply 
connected, irreducible, and has nonnegative 
curvature operator. Moreover, 
$\pi_*(T\R^2)\subset\mathcal D_f$, 
and the normal vector field 
$\delta=n\mathcal H_f/(n-2)$ 
is a Dupin principal normal of 
$f$ satisfying $\pi_*(TN)\subset E_\delta$, 
where $\pi\colon\tilde M^n\to M^n$ 
is the universal covering map. 
\end{enumerate}
In cases $(ii)$--$(iii)$, every compact 
irreducible factor of dimension at least 
three is either a simply connected irreducible 
compact symmetric space, a Riemannian 
manifold with positive curvature operator 
diffeomorphic to a sphere, or a K\"ahler 
manifold biholomorphic to a complex 
projective space. Every compact two-dimensional 
factor is diffeomorphic to $\Sf^2$ 
with nonnegative Gaussian curvature. 
\end{theorem}
\proof
By Proposition \ref{propu}, the manifold $M^n$ 
has nonnegative isotropic curvature. We claim 
that $M^n$ is locally reducible. Otherwise, 
Theorem \ref{acta} would imply that its 
universal cover $\tilde M^n$ is compact, 
contrary to the assumption that 
$\pi_1(M^n)$ is infinite. 

Consequently, Theorem \ref{duke} applies. 
We first show that case $(b)$ of that 
theorem cannot occur. Indeed, suppose 
that $\tilde M^n$ is isometric 
to a Riemannian product 
$(\varSigma^2,g_\varSigma)
\times(N^{n-2},g_N)$, 
where the Gaussian curvature of 
$\varSigma^2$ is negative at some point, 
and $N^{n-2}$ is compact and irreducible. 
The immersion $\tilde f=f\circ\pi$ satisfies 
\eqref{1}. By part $(ii)$ of Proposition
\ref{propprod}, $\tilde M^n$ has nonnegative 
curvature operator. In particular, $\varSigma^2$ 
has nonnegative Gaussian curvature, 
a contradiction. 

Thus, case $(a)$ of Theorem \ref{duke} 
applies, and $\tilde M^n$ is 
isometric to a Riemannian product 
\[
(\R^{n_0},g_0)\times
(N_1^{n_1},g_1)\times\cdots\times
(N_r^{n_r},g_r),
\]
where $n_0\geq1$, $r\geq0$, $g_0$ 
is the Euclidean metric, and each $(N_i,g_i)$ 
is compact, simply connected, 
irreducible Riemannian manifold of 
dimension $n_i\ge2$. 
Moreover, each compact factor of dimension 
two is diffeomorphic to $\Sf^2$ with 
nonnegative Gaussian curvature. 
We distinguish two cases. 
\vspace{1ex}

\noindent\emph{Case I.} 
Suppose that $M^n$ has positive isotropic 
curvature at some point. 

We claim that $n_0=1$ and $r=1$. 
Indeed, if $n_0\geq2$, we may choose, 
at any point of $\tilde M^n$, an 
orthonormal four-frame 
$\{e_i\}_{i=1}^4$ with 
$e_1,e_2\in T\R^{n_0}$. 
The product structure implies that 
the isotropic curvature of this frame 
vanishes, a contradiction. 
Thus, $n_0=1$. Suppose that $r\geq2$. 
Choose an orthonormal four-frame 
$\{e_i\}_{i=1}^4$ such that $e_1\in T\R$, 
$e_2\in TN_1$, and $e_3,e_4\in TN_2$. 
Again, the product structure implies 
that the isotropic curvature of this 
frame vanishes, a contradiction. 

Therefore, $r=1$, and $\tilde M^n$ 
is isometric to a Riemannian 
product $\R\times N^{n-1}$, 
where $N^{n-1}$ is compact, simply 
connected, irreducible, and has 
nonnegative isotropic curvature. 

Moreover, since $M^n$ has positive 
isotropic curvature at some point, 
so does $N^{n-1}$. By Remark $(iv)$ 
in \cite{ses}, the manifold $N^{n-1}$ 
admits a metric of positive isotropic 
curvature. The theorem of Micallef 
and Moore \cite{MM} then implies that 
$N^{n-1}$ is homeomorphic to 
$\Sf^{n-1}$. 
Finally, Theorem \ref{acta} implies 
that $N^{n-1}$ is diffeomorphic 
to $\Sf^{n-1}$. Thus, alternative 
$(i)$ holds. 
\vspace{1ex}

\noindent\emph{Case II.} 
Suppose now that $M^n$ has no points 
of positive isotropic curvature. 
By part $(ii)$ of Proposition 
\ref{propu}, equality holds in 
\eqref{1} at every point. 
\vspace{0.5ex}

\noindent\emph{Subcase $II_a$.} 
Suppose first that $r=0$, so that 
$\tilde M^n$ is isometric to $\R^n$ 
and $M^n$ is flat. Since the ambient 
curvature operator 
is bounded below by $c$, the 
Gauss equation gives 
\[
0=\tau
\geq n(n-1)c+n^2H^2-S.
\]
Consequently, 
\[
n^2H^2\leq S-n(n-1)c.
\]
This together with inequality 
\eqref{1}, which now holds as an equality, 
implies that $c=0$ and $S=0$. Thus, 
$f$ is totally geodesic, and alternative 
$(ii)$ holds. 
\vspace{0.5ex}

\noindent\emph{Subcase $II_b$.} 
Suppose now that $r\geq1$ and 
$n_i<n-2$ for all $0\leq i\leq r$. 
Since the immersion $\tilde f=f\circ\pi$ 
satisfies \eqref{1}, part $(i)$ of Proposition 
\ref{propprod} implies that $\tilde M^n$ has 
nonnegative curvature operator, $c=0$ and 
that $\tilde f$ is totally geodesic. Hence, $f$ 
is also totally geodesic. Moreover, $M^n$ 
has nonnegative curvature operator. By 
Theorem 5 in \cite{chow} and 
the main result in \cite{bohm}, we 
conclude that each irreducible manifold 
$(N_i,g_i)$, $1\leq i\leq r$, 
is a compact symmetric space, 
a Riemannian manifold with positive curvature operator 
diffeomorphic to a sphere, or a K\"ahler manifold 
biholomorphic to a complex projective space. Thus, 
alternative $(ii)$ holds. 
\vspace{0.5ex}

\noindent\emph{Subcase $II_c$.} 
Suppose now that $r\geq1$ and 
$n_i\geq n-2$ for some $0\leq i\leq r$. 
Since every compact factor has 
dimension at least two, the 
only possibilities are 
$$
(n_0,n_1)=(n-2,2),\quad(n_0,n_1)=(2,n-2),
\quad\text{or}\quad
(n_0,n_1)=(1,n-1).
$$

Suppose first that $\tilde M^n$ 
is isometric to a Riemannian 
product $\R^{n-2}\times N_1^2$. 
In this case, part $(ii)$ of Proposition 
\ref{propprod} applies to the immersion 
$\tilde f=f\circ\pi$. Hence, $M^n$ has 
nonnegative curvature operator. 
Moreover, $c=0$, 
$TN_1^2\subset\mathcal D_{\tilde f}$, 
and the normal vector field 
$\tilde\delta=n\mathcal H_{\tilde f}/(n-2)$ 
is a Dupin principal normal of $\tilde f$ 
satisfying $T\R^{n-2}\subset E_{\tilde\delta}$. 
Choosing orthonormal vectors 
$v_1,v_2\in T\R^{n-2}$, 
the Gauss equation gives 
$$
0=K_{\tilde M}(v_1\wedge v_2)
=K_{\bar M}
(\tilde f_*v_1\wedge\tilde f_*v_2)
+\|\tilde\delta\|^2\geq\|\tilde\delta\|^2.
$$
Thus, $\tilde\delta=0$. Consequently, $\tilde f$ 
and thus $f$ are totally geodesic, and 
alternative $(ii)$ holds. 

Suppose now that $\tilde M^n$ is isometric 
to a Riemannian product $\R^2\times N_1^{n-2}$. 
In this case, part $(ii)$ of Proposition 
\ref{propprod} applies to the immersion 
$\tilde f=f\circ\pi$. Hence, $M^n$ has 
nonnegative curvature 
operator. Moreover, $c=0$, 
$T\R^2\subset\mathcal D_{\tilde f}$, 
and the normal vector field 
$\tilde\delta=n\mathcal H_{\tilde f}/(n-2)$ 
is a Dupin principal normal of $\tilde f$ 
satisfying $TN_1\subset E_{\tilde\delta}$. 
Since $\pi$ is a local isometry, 
we obtain 
$\pi_*(T\R^2)\subset\mathcal D_f$ and 
$\pi_*(TN_1)\subset E_\delta$, 
where $\delta=n\mathcal H_f/(n-2)$. 
Thus, alternative $(iii)$ holds. 

Finally, suppose that $\tilde M^n$ is 
isometric to a Riemannian product 
$\R\times N_1^{n-1}$, where $N_1$ 
is a compact, irreducible, simply connected 
manifold. Fix an arbitrary point 
$x_0=(t_0,y_0)\in\tilde M^n$ 
and consider the curve 
$\gamma_0=\tilde f\circ\sigma_0$ and the 
isometric immersion 
$
f_1=\tilde f\circ\sigma_1\colon N_1\to\bar {M}^{n+m},
$ 
where $\sigma_0\colon\R\to\tilde M^n$ and 
$\sigma_1\colon N_1\to\tilde M^n$ 
denote the totally geodesic inclusions 
defined by 
$$
\sigma_0(t)=(t,y_0)\quad \text{and}
\quad \sigma_1(y)=(t_0,y),\quad y\in N_1.
$$

Then Lemma \ref{lemmaprod} applies to 
the immersion $\tilde f$, and \eqref{pro} 
simplifies to 
\begin{align*}%\label{pro1}
\frac{(n-3)}{n-2}\kappa_1^2+
\frac{2(n-3)}{n-2}\sum_{\ell=1}^{n-1}\Big\|\a_{\tilde f}
\big(\sigma_0{_*}\frac{\partial}{\partial t},\sigma_1{_*}e_\ell\big)\Big\|^2\leq
\frac{2(n-3)}{n-2}c+\frac{(n-1)^2}{n-2}H_1^2-S_1\nonumber,
\end{align*} 
where $\kappa_1$ is the first Frenet curvature of 
the curve $\gamma_0$, and $\{e_\ell\}_{\ell=1}^{n-1}$ is an 
orthonormal frame of $TN_1$. 
Since $c\geq0$, the above inequality yields 
\begin{align}\label{pro3}
\frac{(n-3)}{n-2}\kappa_1^2+
\frac{2(n-3)}{n-2}\sum_{\ell=1}^{n-1}\Big\|\a_{\tilde f}
\big(\sigma_0{_*}\frac{\partial}{\partial t},\sigma_1{_*}e_\ell\big)\Big\|^2&\leq4c
+\frac{(n-1)^2}{n-2}H_1^2-S_1\nonumber\\
&\leq4c
+\frac{(n-1)^2}{n-3}H_1^2-S_1. 
\end{align} 

It follows from \eqref{pro3} that the immersion 
$f_1\colon N_1\to\bar {M}^{n+m}$ 
satisfies the pinching condition \eqref{1}. 
Suppose first that $n\geq6$. 
Since $N_1$ is compact, simply connected, 
and irreducible, Theorem \ref{p1fn5} 
implies that either $N_1$ is 
diffeomorphic to $\Sf^{n-1}$, or 
it is a simply connected irreducible 
compact symmetric space or a K\"ahler 
manifold biholomorphic to $\CP^{(n-1)/2}$. 
In the first case, alternative $(i)$ holds.
In the remaining cases, $c=0$ and
$f_1$ is totally geodesic. Thus, 
$S_1=H_1=0$, and \eqref{pro3} yields 
\be\label{k_1=0}
\kappa_1=0,\qquad
\alpha_{\tilde f}
\big(\sigma_0{_*}\frac{\partial}{\partial t},\sigma_1{_*}X\big)=0
\quad\text{for all }X\in TN_1.
\ee
Since $x_0$ was arbitrary, it follows 
that $\tilde f$ is totally geodesic. 
Consequently, $f$ is totally geodesic, 
and alternative $(ii)$ holds. 

Suppose now that $n=5$. Then 
Theorem \ref{n=4} applies to the immersion 
$f_1\colon N_1\to\bar {M}^{5+m}$. 
Since $N_1$ is compact, simply connected, 
and irreducible, only cases $(i)$ and 
$(iii)$-$(a)$ of that theorem can occur. 
In the former case, $N_1$ is diffeomorphic 
to $\Sf^4$, and thus alternative $(i)$ in the 
statement of the theorem holds. 

Assume that case $(iii)$-$(a)$ of 
Theorem \ref{n=4} holds. Then $N_1$ 
is a K\"ahler manifold biholomorphic 
to $\CP^2$, and the immersion $f_1$ 
is $(2,0)$-geodesic with respect to 
its complex structure. Moreover, $f_1$ 
satisfies the pinching condition \eqref{1} 
as an equality at every point, so that 
$S_1=4c+8H_1^2$. 
Substituting this identity into 
\eqref{pro3}, we obtain 
\[
\frac{2}{3}\kappa_1^2+
\frac{4}{3}\sum_{\ell=1}^{4}
\Big\|
\alpha_{\tilde f}
\big(
\sigma_0{_*}\frac{\partial}{\partial t},\sigma_1{_*}e_\ell
\big)
\Big\|^2
\leq -\frac{8}{3}H_1^2.
\]
Consequently, 
\[
H_1=\kappa_1=0,
\qquad
\alpha_{\tilde f}
\big(
\sigma_0{_*}\frac{\partial}{\partial t},\sigma_1{_*}X
\big)=0
\quad\text{for all }X\in TN_1.
\]

Since $\tilde M^5$ is isometric to a Riemannian 
product $\R\times N_1$ and the curve 
$\gamma_0$ is a geodesic, the Gauss 
equation gives 
\[
0=
K_{\bar M}
\Big(
\tilde f_*(\sigma_0{_*}\frac{\partial}{\partial t})
\wedge\tilde f_*(\sigma_1{_*}X)
\Big)\geq c
\]
for every unit vector $X\in TN_1$. 
Thus $c=0$, and hence $S_1=0$. 
It follows that $f_1$ is totally 
geodesic. Together with the 
vanishing of $\kappa_1$ and of the mixed 
second fundamental form, this 
implies that $\tilde f$ is totally 
geodesic. Consequently, $f$ is 
totally geodesic, and alternative 
$(ii)$ holds. 
\vspace{0.5ex}

In alternatives $(ii)$--$(iii)$, 
the classification of the compact 
irreducible factors of dimension 
at least three follows from 
Theorem 5 of \cite{chow}, together 
with the theorem of B\"ohm and 
Wilking \cite{bohm}. Each compact 
two-dimensional factor is diffeomorphic 
to $\Sf^2$ with nonnegative 
Gaussian curvature. 
This completes the proof. 
\qed
\vspace{1ex}

\noindent\emph{Proof of Theorem \ref{n5}:} 
The result follows directly from Theorems \ref{ceven}, 
\ref{p1fn5} and \ref{p1ifn5}.
\qed 

\section{Proof of Theorem \ref{k=1} and submanifolds 
of space forms}

The following lemma, which establishes the relationship 
between the inequality \eqref{2} and the sectional curvature, 
is needed for the proof of Theorem \ref{k=1}.

\begin{lemma}\label{alg22}
Let $f\colon M^n\to\bar M^{n+m}$, $n\geq4$, 
be an isometric immersion, and let $c$ be 
a lower bound on the curvature operator of 
$\bar M^{n+m}$. For any point $x\in M^n$, 
any orthonormal basis $\{e_i\}_{i=1}^n$ 
of $T_xM$, and any orthonormal basis 
$\{\xi_\alpha\}_{\alpha=1}^m$ of $N_fM(x)$, 
the following assertions hold: 
\begin{enumerate}[
label=(\roman*),
topsep=1pt,
itemsep=0.2ex,
leftmargin=*,
partopsep=0.2ex,
parsep=0.2ex
]
\item The sectional curvature satisfies 
\begin{align}\label{ineqK}
K_M(e_1\wedge e_2)
\geq \frac{1}{2}
\Big(2c+\frac{n^2H^2}{n-1}-S\Big)
+\sum_{\alpha=1}^m
\sum_{\substack{i<j\\(i,j)\neq(1,2)}}
(h^\alpha_{ij})^2,
\end{align}
where
$h^\alpha_{ij}=\langle A_{\xi_\alpha}e_i,e_j\rangle$ 
and $A_{\xi_\alpha}$ denotes the shape operator 
associated with $\xi_\alpha$. 

\item If equality holds in \eqref{ineqK}, then 
\begin{equation}\label{ineqalpha}
K_{\bar M}(f_*e_1\wedge f_*e_2)=c,
\qquad
\alpha_{11}+\alpha_{22}
=\alpha_{33}=\cdots=\alpha_{nn},
\end{equation}
where, for convenience, we write 
$\alpha_{ij}=\alpha_f(e_i,e_j)$. 
\end{enumerate}
\end{lemma}

\begin{proof}
We first claim that
\begin{align}\label{ineqK1}
\sum_{\alpha=1}^m
\left(h^\alpha_{11}h^\alpha_{22}
-(h^\alpha_{12})^2\right)
\geq
\frac{1}{2}
\Big(\frac{n^2H^2}{n-1}-S\Big)
+\sum_{\alpha=1}^m
\sum_{\substack{i<j\\(i,j)\neq(1,2)}}
(h^\alpha_{ij})^2.
\end{align}
Moreover, equality implies that
\begin{equation}\label{ineqalpha1}
h^\alpha_{11}+h^\alpha_{22}
=h^\alpha_{33}=\cdots=h^\alpha_{nn}
\end{equation}
for every $1\leq\alpha\leq m$.

Indeed, by the Cauchy--Schwarz inequality, 
\begin{align}\label{min2}
\Big(\sum_{i=1}^n h^\alpha_{ii}\Big)^2
\leq (n-1)
\Big((h^\alpha_{11}+h^\alpha_{22})^2
+\sum_{i=3}^n(h^\alpha_{ii})^2\Big)
=(n-1)\Big(
\sum_{i=1}^n(h^\alpha_{ii})^2
+2h^\alpha_{11}h^\alpha_{22}\Big).
\end{align}
Equivalently,
\begin{equation}\label{min3}
2h^\alpha_{11}h^\alpha_{22}
\geq
\sum_{i\neq j}(h^\alpha_{ij})^2
+\frac{1}{n-1}
(\operatorname{tr}A_{\xi_\alpha})^2
-\|A_{\xi_\alpha}\|^2.
\end{equation}
Summing over $\alpha$ and subtracting 
$\sum_\alpha(h^\alpha_{12})^2$, we obtain 
\begin{align}\label{min4}
\sum_{\alpha=1}^m
\left(h^\alpha_{11}h^\alpha_{22}
-(h^\alpha_{12})^2\right)&\geq
\frac{1}{2}\sum_{\alpha=1}^m
\sum_{i\neq j}(h^\alpha_{ij})^2
+\frac{1}{2(n-1)}
\sum_{\alpha=1}^m
(\operatorname{tr}A_{\xi_\alpha})^2
\nonumber\\
&\quad
-\frac{1}{2}\sum_{\alpha=1}^m
\|A_{\xi_\alpha}\|^2
-\sum_{\alpha=1}^m(h^\alpha_{12})^2
\nonumber\\
&\quad=
\frac{1}{2}
\Big(\frac{n^2H^2}{n-1}-S\Big)
+\sum_{\alpha=1}^m
\sum_{\substack{i<j\\(i,j)\neq(1,2)}}
(h^\alpha_{ij})^2.
\end{align}
This proves \eqref{ineqK1}. Moreover, 
equality in \eqref{ineqK1} forces equality 
in \eqref{min2} for every $1\leq\alpha\leq m$. 
The equality condition in the 
Cauchy--Schwarz inequality then yields 
\eqref{ineqalpha1}. 

On the other hand, the Gauss equation gives 
\begin{align}\label{min44}
K_M(e_1\wedge e_2)
&=K_{\bar M}(f_*e_1\wedge f_*e_2)
+\sum_{\alpha=1}^m
\left(h^\alpha_{11}h^\alpha_{22}
-(h^\alpha_{12})^2\right)
\nonumber\\
&\geq c+
\sum_{\alpha=1}^m
\left(h^\alpha_{11}h^\alpha_{22}
-(h^\alpha_{12})^2\right).
\end{align}
Combining \eqref{min44} with
\eqref{ineqK1} proves \eqref{ineqK}.

Finally, if equality holds in \eqref{ineqK}, 
then equality must hold in both 
\eqref{ineqK1} and \eqref{min44}. 
The former yields the second condition 
in \eqref{ineqalpha}, while the latter 
gives the first. 
\end{proof}
\vspace{1.0ex}

\noindent\emph{Proof of Theorem \ref{k=1}:} 
We first observe that
\begin{equation}\label{12}
S\leq 2c+\frac{n^2H^2}{n-1}
\leq 4c+\frac{n^2H^2}{n-2}
\end{equation}
at every point. Thus \eqref{1} is also 
satisfied, and the conclusions of 
Theorems \ref{n5} and \ref{n=4} apply. 

By Lemma \ref{alg22}, at every point 
$x\in M^n$ and for every orthonormal 
basis $\{e_i\}_{i=1}^n$ of $T_xM$, we have 
\begin{align}\label{ineqKp}
K_M(e_1\wedge e_k)
&\geq\frac{1}{2}
\Big(2c+\frac{n^2H^2}{n-1}-S\Big)+
\sum_{\alpha=1}^m
\sum_{\substack{i<j\\(i,j)\neq(1,k)}}
(h^\alpha_{ij})^2
\end{align}
for every $2\leq k\leq n$. 
In particular, $M^n$ has nonnegative 
sectional curvature. 

Moreover, choosing $e_1=X$ for any 
unit vector $X\in T_xM$ and summing 
\eqref{ineqKp} over $k$, we obtain 
\begin{align}\label{ric}
\Ric_M(X)
&\geq\frac{n-1}{2}
\Big(2c+\frac{n^2H^2}{n-1}-S\Big)
+
\sum_{\alpha=1}^m\sum_{k=2}^n
\sum_{\substack{i<j\\(i,j)\neq(1,k)}}
(h^\alpha_{ij})^2
\nonumber\\
&\geq\frac{n-1}{2}
\Big(2c+\frac{n^2H^2}{n-1}-S\Big)
\geq0.
\end{align}
Hence $M^n$ has nonnegative Ricci 
curvature. We distinguish two cases. 
\vspace{1ex}

\noindent\emph{Case I.}
Suppose that $M^n$ has positive Ricci 
curvature at some point. By the 
Cheeger--Gromoll splitting theorem 
\cite{cg}, the universal cover of $M^n$, 
equipped with the lifted metric, is 
isometric to a Riemannian product 
$\R^p\times N$, where $N$ is compact. 
Since the Ricci curvature is positive 
at some point, we must have $p=0$. 
Thus, the universal cover is compact, 
and hence $\pi_1(M^n)$ is finite. 

Consequently, the alternatives in 
Theorems \ref{n5} and \ref{n=4} 
involving a nontrivial Euclidean factor 
are excluded. If the spherical alternative 
in either theorem holds, then $M^n$ is 
diffeomorphic to a spherical space form, 
or to a sphere if $n=4$. 

In every remaining alternative, equality holds 
in \eqref{1} at every point. Combining this 
with \eqref{12}, we obtain $c=H=0$, 
which in turn implies $S=0$. 
Consequently, $f$ is totally geodesic. 
The Gauss equation then implies that 
the curvature operator of $M^n$ is 
nonnegative. The asserted decomposition 
of its universal cover follows from 
Theorem 5 in \cite{chow} and the main 
result in \cite{bohm}. Thus alternative 
$(ii)$ holds. 
\vspace{1ex}

\noindent\emph{Case II.}
Suppose that at every point $x\in M^n$ 
there exists a unit vector $X\in T_xM$ 
such that $\Ric_M(X)=0$. 

We claim that if, at a point $x$, a unit 
vector $X\in T_xM$ satisfies $\Ric_M(X)=0$, 
then $c=0$, equality holds in \eqref{2} at $x$, 
$X\in\mathcal D_f(x)$, and $\delta=n\mathcal H_f/(n-1)$ 
is a Dupin principal normal of multiplicity at least 
$n-1$ satisfying $E^\perp_\delta(x)\subset\mathcal D_f(x)$. 

To prove the claim, choose an orthonormal 
basis $\{e_i\}_{i=1}^n$ of $T_xM$ with $e_1=X$. 
Since all terms on the right-hand side
of \eqref{ric} are nonnegative, we obtain 
that equality holds in \eqref{2} at $x$, and 
\begin{equation}\label{k1off}
\sum_{\alpha=1}^m\sum_{k=2}^n
\sum_{\substack{i<j\\(i,j)\neq(1,k)}}
(h^\alpha_{ij})^2=0.
\end{equation}
In particular, equality holds in 
\eqref{ineqKp} for every $2\leq k\leq n$. 
Since the sectional curvature is 
nonnegative and $\Ric_M(e_1)=0$, we 
also have 
\[
K_M(e_1\wedge e_k)=0,
\qquad 2\leq k\leq n.
\]

By part $(ii)$ of Lemma \ref{alg22}, 
it follows that 
\begin{equation}\label{ineqalphaki}
K_{\bar M}(f_*e_1\wedge f_*e_k)=c,
\qquad
\alpha_{11}+\alpha_{kk}=\alpha_{ll}
\end{equation}
for all $2\leq k\leq n$ and 
$l\notin\{1,k\}$, where 
$\alpha_{ij}=\alpha_f(e_i,e_j)$. 
Since $n\geq4$, these relations imply 
\begin{equation}\label{k1diag} 
\alpha_{11}=0,\qquad
\alpha_{22}=\cdots=\alpha_{nn}=\delta
\end{equation}
for some normal vector $\delta$. 
On the other hand, \eqref{k1off} implies
that
\begin{equation}\label{k1off2}
\alpha_{ij}=0,\qquad i\neq j.
\end{equation} 

Combining \eqref{k1diag} and 
\eqref{k1off2}, we obtain 
\[
\alpha_{1j}=0,\qquad 1\leq j\leq n,
\]
and
\[
\alpha_{ij}=\delta_{ij}\delta,
\qquad 2\leq i,j\leq n.
\]
Hence, $X\in\mathcal D_f(x)$, and 
$\delta$ is a Dupin principal normal 
of multiplicity at least $n-1$ satisfying 
$E^\perp_\delta(x)\subset\mathcal D_f(x)$. 
In particular, $n\mathcal H_f=(n-1)\delta$. 

It remains to prove that $c=0$. In fact, the 
Gauss equation now gives 
\[
0=K_M(e_1\wedge e_k)
=K_{\bar M}(f_*e_1\wedge f_*e_k)=c.
\]
Thus $c=0$, and this proves the claim. 

If $\delta=0$ everywhere, then 
$S=0$, and hence $f$ is totally 
geodesic. As in \emph{Case I}, the Gauss 
equation and the results of 
\cite{chow,bohm} yield alternative 
$(ii)$. 

Suppose now that $f$ is not totally geodesic. 
Let $\pi\colon\tilde M^n\to M^n$ denote the 
universal covering map. According to the 
Cheeger--Gromoll splitting theorem \cite{cg}, 
the universal cover $\tilde M^n$, 
equipped with the covering metric, splits 
isometrically as a Riemannian product 
$\R^p\times N$, where $N$ is a compact 
manifold. 

We claim that $p\leq1$. Indeed, let 
$x_0\in M^n$ be a point where $f$ is 
not totally geodesic. Then 
$\delta(x_0)\neq0$, and the preceding 
claim gives $\dim\mathcal D_f(x_0)=1$. 
Since every tangent vector to the 
Euclidean factor belongs to the 
relative nullity subspace, we obtain 
$p\leq1$. 

If $p=0$, then $\pi_1(M^n)$ is finite. 
By Theorems \ref{n5} and \ref{n=4}, 
either alternative $(i)$ holds or 
equality holds in \eqref{1} everywhere. 
In the latter case, arguing as in 
\emph{Case I}, we conclude that $f$ is 
totally geodesic, and this is a contradiction. 

Now suppose that $p=1$. Then the universal 
cover of $M^n$ is isometric to a Riemannian 
product $\R\times N^{n-1}$, where $N$ 
is compact and simply connected. 
Moreover, 
$\pi_*(T\R)\subset\mathcal D_f$, and 
$\pi_*(TN)=E_\delta$ at every point 
where $\delta\neq0$. At points where 
$\delta=0$, the immersion $f$ is 
totally geodesic. By the Gauss equation 
and the nonnegativity of the ambient 
curvature operator, we have
\[
\mathcal R_N(x)\geq
\|\delta(\pi(t,x))\|^2 \,{\rm Id}
\]
for every $(t,x)\in\R\times N$. 
In particular, the curvature operator 
of $N$ is nonnegative everywhere. 
Since $f$ is not totally geodesic, 
$\delta$ is nonzero at some point, 
and hence $\mathcal R_N$ is positive 
at some point.

We claim that $N$ 
admits a Riemannian metric with positive 
curvature operator. 
Indeed, let $g(t)$ be the solution of the 
Ricci flow with initial metric $g_N$. 
The nonnegativity of the curvature operator 
is preserved along the flow. Moreover, 
by Hamilton's strong maximum principle 
\cite[Lemma 8.2 and Theorem 8.3]{H}, 
for all sufficiently small $t>0$, the 
kernel of the curvature operator of $g(t)$ 
has constant dimension and is invariant 
under parallel translation. Since the 
curvature operator of $g_N$ is positive 
at some point, the same holds 
for $g(t)$ at that point for sufficiently 
small $t>0$. Consequently, the kernel 
is trivial everywhere, and hence 
$g(t)$ has positive curvature operator. 

Thus $N$ admits a metric with positive 
curvature operator. By the theorem of 
B\"ohm and Wilking \cite{bohm}, $N$ is 
diffeomorphic to a spherical space form. 
Since $N$ is simply connected, it follows 
that $N$ is diffeomorphic to $\Sf^{n-1}$. 
Thus alternative $(iii)$ holds, completing 
the proof.
\qed
\vspace{1ex}

\noindent\emph{Proof of Corollary \ref{sfhigh}:} 
The proof follows immediately from 
Theorem~\ref{n5} and the fact that 
every compact submanifold of Euclidean 
space has a point at which its relative 
nullity subspace is trivial 
(see Corollary~1.6 in \cite{DT}). 
\qed 
\vspace{1ex}

\noindent\emph{Proof of Corollary \ref{cor}:} 
Since Euclidean space admits no compact 
totally geodesic submanifolds, 
alternative~$(ii)$ of Theorem~\ref{k=1} 
is excluded. Alternative~$(iii)$ is also 
impossible, since every compact submanifold 
of Euclidean space has a point at which 
its relative nullity subspace is trivial, 
whereas in alternative~$(iii)$ this subspace 
is nontrivial everywhere. 
Consequently, Theorem~\ref{k=1} 
implies that $M^n$ is diffeomorphic 
to a spherical space form, and to 
$\Sf^4$ when $n=4$. 

Suppose henceforth that $n\geq5$. 
We also have 
\begin{equation}\label{v3}
S\leq 2c+\frac{n^2H^2}{n-1}
\leq nc+\frac{n^3H^2}{2(n-1)}
-\frac{n(n-2)}{2(n-1)}H
\sqrt{n^2H^2+4c(n-1)}.
\end{equation}
Thus, the pinching condition considered 
in \cite{v3} for $k=1$ is satisfied. 
By Theorem 2 of \cite{v3}, either 
$M^n$ is homeomorphic to $\Sf^n$, 
or equality holds throughout \eqref{v3} 
and the submanifold belongs to one 
of the exceptional classes described 
in parts $(i)$--$(iii)$ of that theorem. 
In the latter case, equality throughout
\eqref{v3} forces $c=0$. However, none 
of the exceptional cases described in 
parts $(i)$--$(iii)$ of Theorem~2 in 
\cite{v3} is compatible with our 
pinching condition when $c=0$. 
Consequently, $M^n$ is homeomorphic 
to $\Sf^n$ and hence simply connected. 
Since $M^n$ is also diffeomorphic to 
a spherical space form, it follows that 
$M^n$ is diffeomorphic to $\Sf^n$. 
\qed

\section{Examples of submanifolds satisfying the pinching condition}

In this section, we construct several 
families of submanifolds satisfying 
\eqref{1}, illustrating the range of 
geometric and topological configurations 
allowed by our results. In particular, 
we exhibit geometrically distinct 
immersions of manifolds diffeomorphic 
to $\Sf^n$ or $\Sf^1\times\Sf^{n-1}$, 
including examples satisfying the 
pinching condition strictly and others 
realizing equality. These examples 
also show that our rigidity conclusions 
cannot, in general, be strengthened 
to extrinsic rigidity. 

There are many compact hypersurfaces 
in Euclidean space that satisfy \eqref{1} 
and are diffeomorphic to spheres. 
To illustrate this, let 
$f\colon M^n\to\R^{n+1}$, $n\geq4$, 
be an ovaloid with principal curvatures 
$0<\lambda_1\leq\cdots\leq\lambda_n$. 
It follows from the inequalities 
$S\leq n\lambda^2_n$ and $H\geq\lambda_1$ 
that, if the global minimum of $\lambda_1$ and 
the global maximum of $\lambda_n$ satisfy 
\be\label{lambda2}
\max_M\lambda _n\leq\min_M\lambda _1\big({\frac{n}{n-2}}\big)^{1/2}, 
\ee
then $f$ satisfies condition \eqref{1} at every point. 
In particular, if \eqref{lambda2} holds strictly, then 
\eqref{1} is satisfied as a strict inequality everywhere. 
By Hadamard’s classical theorem, every ovaloid in 
$\R^{n+1}$ is diffeomorphic to $\Sf^n$. 

A large class of ellipsoids provides concrete 
examples satisfying condition \eqref{lambda2}. 
Consider, for instance, the 
ellipsoid in $\R^{n+1}$ defined by 
$$ 
\frac{x_1^2}{a_1^2}+\cdots+
\frac{x_{n+1}^2}{a_{n+1}^2} = 1,
$$
where $0<a_1\leq\dots\leq a_{n+1}$. A 
straightforward computation shows that the 
minimum and the maximum of the principal 
curvatures of the ellipsoid are $a_1/a_{n+1}^2$ 
and $a_{n+1}/a_1^2$, respectively. It follows 
that condition \eqref{lambda2} holds provided 
$$
a_{n+1}\leq a_1\big({\frac{n}{n-2}}\big)^{1/6}. 
$$

The next proposition presents geometrically 
distinct isometric immersions of manifolds diffeomorphic 
to the product $\Sf^1 \times \Sf^{n-1}$, providing further 
examples that satisfy condition \eqref{1}. 

\begin{proposition}\label{ex}
Let $g\colon N^{n-1}\to\R^{m_1}$,
$n\geq4$, be an isometric immersion
of a compact Riemannian manifold satisfying
\[
S_g<\frac{(n-1)^2H_g^2}{n-2}
\]
at every point. Let 
$\gamma\colon\Sf^1\to\R^{m_2}$ 
be a regular closed curve, whose 
first curvature $\kappa_1$ satisfies 
\begin{equation}\label{k1}
\max_{\Sf^1}\kappa_1^2\leq
\frac{n-2}{n-3}
\min_N\left\{
\frac{(n-1)^2H_g^2}{n-2}-S_g
\right\}.
\end{equation}
Then the product immersion 
\[
f=\gamma\times g\colon
\Sf^1\times N^{n-1}
\to\R^{m_1+m_2}
\]
satisfies \eqref{1}. 
Moreover, if $n\geq5$, then 
$N^{n-1}$ is diffeomorphic 
to $\Sf^{n-1}$. 
\end{proposition}
\begin{proof}
The squared norm of the second 
fundamental form and the squared 
mean curvature of $f$ satisfy 
\[
S_f=\kappa_1^2+S_g,
\qquad
n^2H_f^2=\kappa_1^2
+(n-1)^2H_g^2.
\]
Consequently, 
\[
\frac{n^2H_f^2}{n-2}-S_f
=
\frac{(n-1)^2H_g^2}{n-2}-S_g
-\frac{n-3}{n-2}\kappa_1^2.
\]
Thus \eqref{k1} is equivalent to 
the pinching condition \eqref{1} 
for the product immersion. 

Finally, for $n\geq5$, 
Corollary \ref{cor} applied to $g$ 
implies that $N^{n-1}$ is 
diffeomorphic to $\Sf^{n-1}$. 
\end{proof}
\vspace{1ex}

We next construct a large class of 
immersions satisfying the hypotheses 
of Proposition \ref{ex}. These include 
geometrically distinct ovaloids in 
Euclidean space, which in turn yield 
families of product immersions 
satisfying \eqref{1}, either strictly 
or with equality. 

Let $g\colon N^{n-1}\to\R^n$, $n\geq4$, 
be an ovaloid with principal curvatures 
$0<\lambda_1\leq\dots\leq\lambda_{n-1}$. 
Since $S_g\leq (n-1)\lambda^2_{n-1}$ and 
$H_g\geq\lambda_1$, it follows that 
$g$ satisfies the hypotheses of 
Proposition \ref{ex} provided that 
\be\label{lambda}
\max_N\lambda _{n-1}<
 \min_N\lambda _1\big({\frac{n-1}{n-2}}\big)^{1/2}.
\ee 

A large class of ellipsoids provides concrete 
examples satisfying condition \eqref{lambda}. 
Specifically, the ellipsoid in $\R^n$ defined by 
$$ 
\frac{x_1^2}{a_1^2}+\cdots+
\frac{x_n^2}{a_n^2} = 1,
$$
where $0<a_1\leq\dots\leq a_n$, satisfies 
condition \eqref{lambda} provided that 
$$
a_n<a_1\big({\frac{n-1}{n-2}}\big)^{1/6}. 
$$

Consequently, Proposition \ref{ex} 
yields numerous geometrically distinct 
immersions of manifolds diffeomorphic 
to $\Sf^1\times\Sf^{n-1}$ 
satisfying \eqref{1}. By choosing 
$\gamma$ so that the inequality 
\eqref{k1} is strict, we obtain 
compact examples satisfying \eqref{1} 
strictly at every point. Moreover, 
the strict pinching condition is 
preserved under sufficiently small 
smooth deformations of these immersions. 

\begin{example}[\bf{Totally geodesic submanifolds of 
symmetric spaces of compact type}] 
\emph{Every compact totally geodesic 
submanifold of a symmetric space 
of compact type satisfies \eqref{1} 
with $c=0$. The classification of 
totally geodesic submanifolds in 
symmetric spaces is a classical 
problem in submanifold geometry, 
which remains open in full generality. 
Complete classifications are known 
for irreducible symmetric spaces of rank 
one and two, through the work of Wolf 
\cite{Wolf}, Chen and Nagano \cite{CN1,CN2}, 
and Klein \cite{Klein1,Klein2,Klein3,Klein4}. 
The classification has also been established 
for products of rank-one symmetric spaces 
\cite{RV}. More recently, Kollross and 
Rodr\'iguez-V\'azquez \cite{KRV} classified 
maximal totally geodesic submanifolds 
of exceptional symmetric spaces. 
Nevertheless, the general classification 
problem for symmetric spaces of compact 
type of higher rank remains open. 
}
\end{example}

\begin{example}[\bf{Manifolds with 
nonnegative curvature operator}]
\emph{
Let $M^n$ be a compact Riemannian manifold with 
nonnegative curvature operator, and let 
$\bar M^{n+m}=M^n\times N^m$ be the Riemannian 
product of $M^n$ with an arbitrary Riemannian 
manifold $N^m$ whose curvature operator is also 
nonnegative. For any fixed point $x_0\in N^m$, 
consider the totally geodesic inclusion 
$$
f\colon M^n\to\bar {M}^{n+m},\quad f(x)=(x, x_0).
$$ 
Since the second fundamental form of $f$ 
vanishes identically, condition~\eqref{1} 
is satisfied with equality for $c=0$. 
}
\end{example}

\begin{example}\emph{
Consider the product immersion 
\[
f=\mathrm{Id}_{\mathbb T^2}\times j
\colon \mathbb T^2\times\Sf^{n-2}(r)
\to\mathbb T^2\times\R^{n-1},
\]
where $\mathbb T^2=\Sf^1\times\Sf^1$ 
is a flat $2$-torus, 
$\mathrm{Id}_{\mathbb T^2}$ denotes 
its identity map, and 
$j\colon\Sf^{n-2}(r)\to\R^{n-1}$ 
is the standard umbilical inclusion. 
A direct computation shows that 
the pinching condition~\eqref{1} 
holds with equality for $c=0$. 
This immersion realizes the configurations 
described in alternative~$(iii)$-$(b)$ 
of Theorem~\ref{n5} for $n\geq5$ and 
alternative~$(iii)$-$(c)$ of 
Theorem~\ref{n=4} for $n=4$. 
}

\emph{ The standard inclusion
$\Sf^1(r)\times\Sf^{n-1}(r)
\hookrightarrow\Sf^1(r)\times\R^n$ 
provides an example realizing
alternative $(ii)$-$(b)$ of
Theorem \ref{k=1}.
}
\end{example}

\begin{example}[\bf{Totally geodesic submanifolds of flat manifolds}]\emph{
Every compact, positive-dimensional, 
totally geodesic submanifold of a 
flat Riemannian manifold satisfies 
\eqref{1} with $c=0$ and has infinite 
fundamental group. Recall that every 
compact flat Riemannian manifold is 
isometric to a quotient $\R^{n+m}/\Gamma$, 
where $\Gamma$ is a Bieberbach 
group (see \cite{W}). More precisely, 
every connected, compact, embedded 
totally geodesic submanifold of 
$\R^{n+m}/\Gamma$ lifts to an affine 
subspace $A\subset\R^{n+m}$ whose 
stabilizer 
\[
\Gamma_A
=\{\gamma\in\Gamma:\gamma(A)=A\}
\]
acts cocompactly on $A$. 
The submanifold is then isometric 
to the compact flat manifold 
$A/\Gamma_A$.
}
\end{example}

\begin{example}\emph{
The standard isometric embedding 
$\Sf^2(r)\times\Sf^2(\sqrt{R^2-r^2})
\hookrightarrow 
\Sf^5(R)$ 
satisfies condition \eqref{1} with equality 
for $c=1/R^2$, thereby realizing the 
configuration described in alternative 
$(iii)$-$(b)$ of Theorem~\ref{n=4}. 
}

\emph{
Similarly, the 
standard isometric embedding 
$\Sf^1(r)\times\Sf^3(\sqrt{1-r^2})\hookrightarrow\Sf^5$ 
satisfies \eqref{1} with $c=1$ whenever $r\geq1/2$. 
Since its fundamental group is 
infinite, this example shows 
that alternative $(ii)$ of 
Theorem \ref{n=4} does occur. 
}
\end{example}

\begin{example}[\bf{$(2,0)$-geodesic K\"ahler submanifolds}]
\emph{
The standard minimal embedding 
$g\colon\CP^2\to\Sf^7$ 
of the complex projective plane equipped 
with the Fubini--Study metric of constant 
holomorphic sectional curvature $4/3$ 
is a $(2,0)$-geodesic immersion satisfying 
\eqref{1} with equality for $c=1$, realizing 
alternative $(iii)$-$(a)$ of 
Theorem \ref{n=4}. 
}
\end{example}

\noindent Theodoros Vlachos\\
University of Ioannina \\
Department of Mathematics\\
45110 Ioannina -- Greece\\
e-mail: tvlachos@uoi.gr

\end{document}